\documentclass[oper,nonblindrev]{informs4} % current default for manuscript submission
\OneAndAHalfSpacedXII % Current default line spacing
\usepackage{natbib}
\usepackage{multirow}
\usepackage{booktabs}
\usepackage{placeins}
 \bibpunct[, ]{(}{)}{,}{a}{}{,}%
 \def\bibfont{\small}%
\TheoremsNumberedThrough     % Preferred (Theorem 1, Lemma 1, Theorem 2)
\EquationsNumberedThrough    % Default: (1), (2), ...
\ECRepeatTheorems
\usepackage{amsmath}
\usepackage{graphicx}
\usepackage{bbm}
\usepackage{amssymb} 
\usepackage{amsfonts}
\usepackage{datetime}
\usepackage{mathtools}
\usepackage{comment}
\usepackage{xcolor}
\usepackage{threeparttable}
\usepackage{algorithm}
\usepackage{algorithmic}
\usepackage{caption}
\usepackage{subcaption}
\usepackage{endnotes}
\usepackage{enumitem}
\usepackage{booktabs}

\usepackage[colorinlistoftodos]{todonotes}
\definecolor{links}{RGB}{204,36,29}
\usepackage[colorlinks=true,breaklinks=true,bookmarks=true,urlcolor=links,citecolor=links,linkcolor=links,bookmarksopen=false,draft=false]{hyperref}

\DeclareMathOperator{\Span}{span}

\DeclareMathOperator{\Sign}{sign}

\DeclareMathOperator{\Orth}{orth}

\newcommand{\R}{\mathbb{R}}

\def \p {\boldsymbol{p}}

\def \q {\boldsymbol{q}}
\def \X {\boldsymbol{X}}
\def \Y {\boldsymbol{Y}}
\def \Z {\boldsymbol{Z}}

\def \x {\boldsymbol{x}}
\def \blambda {\boldsymbol{\lambda}}

\def \p {\boldsymbol{p}}
\def \q {\boldsymbol{q}}
\def \A {\mathbf{A}}
\def \B {\mathbf{B}}
\def \bphi {\boldsymbol{\phi}}
\def \x {\boldsymbol{x}}

\def \z {\boldsymbol{z}}
\def \X {\boldsymbol{X}}
\def \Y {\boldsymbol{Y}}
\def \Z {\boldsymbol{Z}}

\DeclareMathOperator{\orth}{orth}
\DeclareMathOperator{\sign}{sign}
\def \blambda {\boldsymbol{\lambda}}

\makeatletter
\newcommand*{\rom}[1]{\expandafter\@slowromancap\romannumeral #1@}
\makeatother

\def \z {\boldsymbol{z}}

\DeclareMathOperator{\conv}{conv}

\DeclareMathOperator{\cone}{cone}

\DeclareMathOperator{\for}{ for }

\DeclareMathOperator{\HM}{\mathsf{HM}} %mathtt

\DeclareMathOperator{\AI}{\mathsf{AI}}

 \renewcommand{\theARTICLETOP}{}

\begin{document}
%%%%%%%%%%%%%%%%

%\setcounter{page}{1}
%\DOI{10.1287/xxxx.0000.0000}%


% \RUNTITLE{}

% \TITLE{Bundling Information Goods of Decreasing Value}
\TITLE{Revisiting the Houtman-Maks Index: Combinatorial Afriat Inequalities and Compact Formulations}

\ARTICLEAUTHORS{%
\AUTHOR{Taotao He}
\AFF{Antai College of Economics and Management, Shanghai Jiao Tong University, \href{mailto:hetaotao@sjtu.edu.cn}{hetaotao@sjtu.edu.cn},~\href{https://taotaoohe.github.io/}{https://taotaoohe.github.io/}}
\AUTHOR{Runfa Hu}
\AFF{Antai College of Economics and Management, Shanghai Jiao Tong University, \href{mailto:rova_hu@sjtu.edu.cn}{rova\_hu@sjtu.edu.cn},~\href{https://runfahu.github.io/}{https://runfahu.github.io/}}
% \AUTHOR{Yating Zhang}
% \AFF{Antai College of Economics and Management, Shanghai Jiao Tong University, \href{mailto:ytzhang20@sjtu.edu.cn}{ytzhang20@sjtu.edu.cn}}
}  

\ABSTRACT{% 
The Houtman--Maks index is a fundamental metric for quantifying violations of the Generalized Axiom of Revealed Preference (GARP) in consumer choice datasets. Because computing this index is NP-hard, exact methods rely on mixed-integer programming (MIP). However, the state-of-the-art formulation scales quadratically in the number of binary and continuous variables relative to the dataset size and depends on weak big-$M$ constants, severely limiting its computational scalability. In this technical note, we derive a big-$M$-free formulation that requires only a linear number of binary and continuous variables. At the core of our approach is a combinatorial version of Afriat inequalities, which yields a potential-based characterization of GARP without auxiliary scaling multipliers. To establish this characterization, we introduce the notion of orthant hull and derive its combinatorial representation via hypersum operations on sign vectors from~\cite{baker2019matroids}. Convexifying pairwise potential inequalities with indicator variables then yields our compact MIP formulation. Computationally, our approach scales to large-scale datasets; theoretically, it extends to general budget sets, finite choice sets, and symmetric rationalizability.
% Numerical experiments show that our formulation dramatically tightens the root LP relaxation, solves substantially more instances to optimality, and achieves an over 40x speedup on perturbed Cobb--Douglas datasets relative to the state-of-the-art benchmark. Finally, our framework extends naturally to general budget sets, finite choice sets, and symmetric rationalizability.

% Individual choice is not only determined by idiosyncratic evaluations of alternatives but also by the influence of the observed behavior of others. However, it is not apparent that assortment planning under peer influence (API) is mixed-integer linear programming (MILP) representable. This paper confirms the MILP representability of API. Moreover, we propose a formulation technique to obtain compact formulations for API. Last, we demonstrate that such technique finds applications in A, B and C.
}

%\KEYWORDS{deterministic inventory theory; infinite linear programming duality; 
%  existence of optimal policies; semi-Markov decision process; cyclic schedule}

\HISTORY{This version: \today}

\maketitle
%%%%%%%%%%%%%%%%%%%%%%%%%%%%%%%%%%%%%%%%%%%%%%%%%%%%%%%%%%%%%%%%%%%%%%

% Samples of sectioning (and labeling) in MSOM
% NOTE: (1) \section and \subsection do NOT end with a period
%       (2) \subsubsection and lower need end punctuation
%       (3) capitalization is as shown (title style).
%
%\section{Introduction.}\label{intro} %%1.
%\subsection{Duality and the Classical EOQ Problem.}\label{class-EOQ} %% 1.1.
%\subsection{Outline.}\label{outline1} %% 1.2.
%\subsubsection{Cyclic Schedules for the General Deterministic SMDP.}
%  \label{cyclic-schedules} %% 1.2.1
%\section{Problem Description.}\label{problemdescription} %% 2.

%%%%%%%%%%%%%%%%%%%%%%%%%%%%%%%%%%%%%%%%%%%%%%%%%
%%%%%%%%%%%% SECTION 1: Introduction %%%%%%%%%%%%
%%%%%%%%%%%%%%%%%%%%%%%%%%%%%%%%%%%%%%%%%%%%%%%%%
\section{Introduction}\label{sec:intro}

Utility maximization is a foundational assumption in economic theory, spanning consumer choice, general equilibrium, and macroeconomic models. Although standard in theoretical analysis, an agent's underlying utility function or preference relation is unobservable in practice. To address this challenge, revealed preference theory, pioneered by \cite{Samuelson1938}, poses an \textit{inverse problem}: given a finite dataset of observed prices and choices, when is the observed behavior consistent with utility maximization, and if so, how can we recover a rationalizing utility function from the data? 

% can we conclude that the agent can be seen as a utility maximizer? In general, this question asks for a rule to test whether the observed choices are consistent with utility maximization.

% While such a rule may  depend on the specific assumptions on the utility function and the choice set, we focus on the classical settings in consumer theory, where the choice set is a linear budget set and the utility function is strictly increasing and concave. We will formalize this consistency question as a rationalizability problem in the next section.

\citet{Afriat1967} shows that the Generalized Axiom of Revealed Preference (GARP) and Afriat's inequalities (AI) each provide a necessary and sufficient condition for rationalizability by a locally nonsatiated utility, with any feasible solution to AI yielding an explicit utility construction. Computationally, \citet{Varian1982} gives a fast combinatorial algorithm for testing GARP, whereas a feasible solution to AI can be computed in polynomial time via linear programming. Subsequent proofs and algorithms have been developed using diverse techniques, including nonnegative cycles~\citep{piaw2003afrait}, linear programming duality~\citep{fostel2004two}, the minimax theorem~\citep{geanakoplos2013afriat}, and the separating hyperplane theorem~\citep{beggs2021afriat}. Revealed preference tests have also been extended to various settings, including weakly separable utility~\citep{cherchye2015revealed}, collective choice~\citep{nobibon2016revealed}, Pareto-optima characterization~\citep{che2024near}, dynamic environments~\citep{chen2020online} and stochastic rationalizability~\citep{smeulders2021nonparametric}. For broader surveys of revealed preference, see~\citet{chambers2016revealed} for theoretical foundations and~\citet{smeulders2019revealed} for algorithmic perspectives.
% Both the theoretical and computational aspects of this problem have been extensively studied by economists and computer scientists in the literature. Theoretically, the non-parametric rule for the classical settings in consumer theory has been well established by Afriat theorem \citep{Afriat1967, Varian1982}, which will be introduced in the next section. Moreover, the problem has been generalized to various settings: the choice set can be non-linear and discrete \citep{forges2009afriat, cosaert2015revealed}, and the utility function can be subject to restrictions such as symmetry \citep{chambers2018note}. We discuss these extensions in Section~\ref{sec:extensions}. Computationally, both efficient algorithms and computational complexity results have been established for various settings \citep{Smeulders2019}. In the classical consumer setting, a polynomial-time test for the consistency question has been established \citep{Varian1982, talla2015note}.

% However, a limitation of such efficient algorithms is that they only provide a binary answer to the consistency question, i.e., whether the observed choices are consistent with utility maximization. 

However, observed choices are often noisy and inconsistent, motivating goodness-of-fit indices that measure departures from utility maximization. In this technical note, we focus on the Houtman--Maks (HM) index \citep{houtman1985determining}, defined as the minimum fraction of observations whose removal leaves the remaining choices consistent with utility maximization. Imposing no structural assumptions on the utility function or choice set, the HM index applies across diverse choice environments and sees wide empirical use \citep{choi2007consistency,nitsch2022reliability,chen2023emergence}. Yet, computing it is NP-hard and hard to approximate \citep{smeulders2014goodness}. Exact mixed-integer programming (MIP) formulations have been developed \citep{heufer2015consistent,demuynck2023computing}; notably, the state-of-the-art formulation of \citet{demuynck2023computing} marks a major practical breakthrough, solving several classic datasets within one second and making exact computation accessible for empirical work. However, existing formulations require quadratically many binary and continuous variables and rely on big-$M$ constants, limiting their scalability to large datasets.

Two primary factors governing the computational performance of an MIP formulation are its size and the strength of its continuous relaxation, obtained by dropping all integrality requirements~\citep{vielma2015mixed}. In this technical note, we derive a strong, big-$M$-free formulation that requires only a linear number of binary and continuous variables.

To derive a compact formulation, we revisit Afriat inequalities, which were initially proposed to recover a utility function using two sets of continuous variables: level variables $\bphi$ and scaling multipliers $\blambda$~\citep{Afriat1967}.  Although Afriat inequalities can test consistency, they carry extra variables for utility recovery. This observation prompts a natural question:
\begin{equation}\label{eq:question}
\text{\itshape Can we eliminate the scaling multipliers $\blambda$ to obtain a leaner consistency test?}
\end{equation}
We address this question by proposing a combinatorial version of Afriat inequalities (Theorem~\ref{thm:afr_ext}), which leads to a potential-based characterization without scaling multipliers (Corollary~\ref{cor:potential}). These characterizations are established by reducing GARP testing to checking disjointness between the negative orthant and the cyclic embedding of the consumption dataset. Our main technical insight is that this disjointness is preserved under the orthant hull operation (Proposition~\ref{prop:nonempty}),  that is, taking the convex hull of the union of orthants generated by points in the cyclic embedding. This is an intriguing property of the cyclic embedding, as such preservation does not hold for general sets.

One of the most successful techniques for deriving strong MIP formulations is convexification~\citep{conforti2014integer}, which typically proceeds by identifying key problem substructures and characterizing their convex hulls. Following this principle, we derive an explicit convex hull description for pairwise potential inequalities coupled with indicator variables (Theorem~\ref{thm:convex-hull}). Combining our potential-based characterization with this convexification result, we obtain a compact, big-$M$-free MIP formulation that achieves substantial computational improvements over the state-of-the-art formulation of \citet{demuynck2023computing} and extends naturally to general budget sets~\citep{forges2009afriat}, finite choice sets~\citep{cosaert2015revealed}, and symmetric rationalizability~\citep{chambers2018note}.

% Convexification is a cornerstone technique for deriving strong MIP formulations~\citep{conforti2014integer}, typically proceeding by identifying key problem substructures and characterizing their convex hulls. Following this principle, we derive an explicit convex hull description for pairwise potential inequalities coupled with indicator variables (Theorem~\ref{thm:convex-hull}), which strengthens the continuous relaxation without introducing big-$M$ parameters. Combining our potential-based characterization with this convexification result, we obtain a compact, big-$M$-free MIP formulation that achieves substantial computational improvemetns over the state-of-the-art formulation of \citet{demuynck2023computing}, and extends naturally to general budget sets, finite choice sets, and symmetric rationalizability.

The remainder of this note is organized as follows. Section~\ref{sec:preliminaries} presents preliminaries on GARP, Afriat inequalities, and the HM index. Section~\ref{sec:characterization} establishes our combinatorial version of Afriat inequalities and the resulting potential-based characterization. Section~\ref{sec:formulations} derives the pairwise convex hull description and our MIP formulation of the HM index. Section~\ref{sec:computational-results} presents preliminary computational results, and we refer the reader to Section~\ref{sec:extensions} for theoretical extensions.

%%%%%%%%%%%%%%%%%%%%%%%%%%%%%%%%%%%%%%%%%%%%%%%%%%%
%%%%%%%%%%%% SECTION 2： Preliminaries %%%%%%%%%%%%
%%%%%%%%%%%%%%%%%%%%%%%%%%%%%%%%%%%$%%%%%%%%%%%%%%%
\section{Preliminaries}\label{sec:preliminaries}

% \subsection{Consumer Decision and Rationalizability}

In the classical consumer setting, observed choices are represented by a finite consumption dataset of $n$ price-quantity pairs $D := \{ (\p^i, \q^i) \}_{i \in [n]}$ across $m$ goods, where for each observation $i$,  the price vector $\p^i \in \R^m_{>}:=\{\p \in \R^m \mid p_j >0 \for j \in [m] \}$ and the consumption bundle $\q^i \in \R^m_{\geq}:= \{\q \in \R^m \mid q_j \geq 0 \for j \in [m]\}$. Given the computation dataset, the problem of whether it is consistent with utility maximization can be formulated as the following rationalizability problem:

% For the classical settings in consumer theory, the observed choices are usually represented by a finite set of price-quantity pairs called a consumption dataset. More specifically, suppose there are $m$ kinds of different goods in the market, and we observe a set of $n$ price-quantity pairs $\{ (\p^i, \q^i) \}_{i \in [n]}$, where the vector $\p^i \in \mathbb{R}^m_{>0}$ records the prices of $m$ goods in the $i$-th observation, and the vector $\q^i \in \mathbb{R}^m_{\ge 0}$ records the quantities of the goods purchased by the agent in the $i$-th observation. Here, $\mathbb{R}^m_{>0}$ and $\mathbb{R}^m_{\ge 0}$ denote the strictly positive and the non-negative orthants of $\mathbb{R}^m$, respectively, meaning that the agent can purchase zero or positive quantities of the goods, but the prices are strictly positive. 

% For the utility function, we focus on the class of strictly increasing and concave utility functions, which is a standard assumption in consumer theory. Specifically, the utility function $u: \mathbb{R}^m_{\ge 0} \to \mathbb{R}$ assigns a real number to each bundle of goods, and it is strictly increasing if, for any $\q,\q'$ in its domain, $u(\q)>u(\q')$ whenever $q_i\ge q'_i$ for all $i$ and $\q\ne\q'$. 

% Under the above settings, the problem of whether the observed choices are consistent with utility maximization can be formulated as the following \textbf{rationalizability problem}:
\begin{definition}[Rationalizability]
    A consumption dataset $D = \{ (\p^i, \q^i) \}_{i \in [n]}$ is said to be rationalizable  if there exists a strictly increasing function utility function $u: \mathbb{R}^m_{\geq} \to \mathbb{R}$ such that, for each observation $i$, 
    \begin{equation*}
        \q^i \in \arg\max_{\q} \Bigl\{ u(\q) \Bigm| \langle \p^i, \q \rangle \le \langle \p^i, \q^i \rangle,\ \q  \geq 0 \Bigr\}.
    \end{equation*}
\end{definition}
This definition is intuitive: if $D$ is rationalizable, each chosen bundle $\q^i$ maximizes utility over the budget set $\{ \q \in \mathbb{R}^m \mid \langle \p^i, \q \rangle \le \langle \p^i, \q^i \rangle, \q \geq 0  \}$ at prices $\p^i$. 

\cite{Afriat1967} characterizes when rationalizability holds via the Generalized Axiom of Revealed Preference (GARP) and Afriat inequalities. For a consumption dataset $D$, define its \textit{revealed preference matrix} $\A := (a_{ij}) \in \mathbb{R}^{n \times n}$ by 
\begin{equation}\label{eq:revealedmatrix}
    a_{ij} := \langle \p^i, \q^j - \q^i \rangle \quad \for i, j \in [n].
\end{equation}
An entry $a_{ij} \le 0$ indicates that $\langle \p^i, \q^j \rangle \le \langle \p^i, \q^i \rangle$, meaning bundle $\q^j$ was affordable when $\q^i$ was chosen. So, $\q^i$ is directly revealed preferred to $\q^j$, and strictly so if $a_{ij} < 0$.

\begin{definition}[GARP]
A consumption dataset $D = \{ (\p^i, \q^i) \}_{i \in [n]}$ satisfies GARP if, for every sequence of indexes $i_1, i_2 , \ldots, i_L$ with $i_{L+1} = i_1$ such that $a_{i_k i_{k+1}} \leq 0$ for $k \in [L]$ 
% \begin{equation*}
%     a_{i_k i_{k+1}} \leq 0 \quad \for k \in [L],
% \end{equation*}
then all of the above inequalities hold with equality.
\end{definition}

The second characterization states that $D$ is rationalizable if and only if there exist level variables $\bphi \in \mathbb{R}^n$ and scaling multipliers $\blambda \in \mathbb{R}^n$ such that the following system is feasible:
\begin{equation*}\label{eq:afriat}
\begin{aligned}
    & \phi_i - \phi_j + \lambda_i a_{ij} \ge 0 \quad \for i,j \in [n] \quad \text{ and } \quad \lambda_i > 0 \quad \for i \in [n].
\end{aligned}
    \tag{\textsc{AI}}
\end{equation*}

\begin{theorem}[Afriat's Theorem]
    Given a consumption dataset $D = \{ (\p^i, \q^i) \}_{i \in [n]}$, the following three statements are equivalent:
    \begin{enumerate}
        \item $D$ is rationalizable;
        \item $D$ satisfies GARP;
        \item The system \eqref{eq:afriat} is feasible.
    \end{enumerate}
\end{theorem}
Our combinatorial Afriat inequalities (Theorem~\ref{thm:afr_ext}) and potential-based characterization (Corollary~\ref{cor:potential}) contribute to the longstanding literature on proofs and algorithms for Afriat's Theorem~\citep{Afriat1967,Varian1982,piaw2003afrait,fostel2004two,geanakoplos2013afriat,beggs2021afriat}. Whereas prior work builds directly on Afriat's system in level variables $\bphi$ and scaling multipliers $\blambda$, we depart from this paradigm by resolving Question~\ref{eq:question}, eliminating the scaling multipliers altogether to yield a purely potential-based characterization.

% In Section \ref{sec:characterization}, we propose combinatorial Afriat inequalities, which subsumes~\eqref{eq:afriat}, and show that such inequalities yields a valid  characterization, generalizing Afriat's Theorem. 

% provide an alternative characterization by using a combinatorial Afrian reprove Afriat theorem and provide a combinatorial extension of \eqref{eq:afriat}, which serves as the foundation of our MIP formulation for computing the HMI. 

% \subsection{The Houtman-Maks Index}
The Houtman-Maks (HM) index measures departure from rationalizability via the largest rationalizable subset of a consumption dataset $D$. Intuitively, more severe violations require more deletions, leaving a smaller rationalizable subset. Formally, the index is defined by
\begin{equation*}\label{prob:HMI}
    \HM(D) := \max_{R\subseteq D}\biggl\{ \frac{|R|}{|D|}  \biggm| R \ \text{is rationalizable} \biggr\},
    \tag{\textsc{HM}}
\end{equation*}
where $\HM(D)=1$ indicates that $D$ is rationalizable. 

Our MIP formulation of $\HM(D)$ contributes to exact MIP formulations of revealed preference tests~\citep{cherchye2015revealed,nobibon2016revealed,heufer2015consistent,demuynck2023computing}. Viewed through the lens of data-driven inverse optimization~\citep{chan2025inverse}, computing the HM index corresponds to nonparametric utility imputation under an $\ell_0$ cardinality loss, identifying the minimum number of outlier observations whose removal renders the remaining choices exactly optimal for a locally nonsatiated utility. Our linear-size, big-$M$-free formulation provides a scalable exact MIP method for this nonparametric inverse optimization problem. 
% A simple but important fact about this index is that it does not ask for an explicit utility function. It only asks whether a subset of data is rationalizable or not, which indicates its combinatorial property. This fact is simple but crucial for our formulation. Therefore, we use $\HM(D)=1$ to describe the rationalizability of $D$. 

% Computationally, \citet{demuynck2023computing} proposed the state-of-the-art MIP formulation for computing the HMI (see Appendix~\ref{app:dr-formulation}).

%%%%%%%%%%%%%%%%%%%%%%%%%%%%%%%%%%%%%%%%%%%%%%%%%%%%%
%%%%%%%%%%%% SECTION 3: Characterization %%%%%%%%%%%%
%%%%%%%%%%%%%%%%%%%%%%%%%%%%%%%%%%%%%%%%%%%%%%%%%%%%%
\section{Combinatorial Afriat inequalities}\label{sec:characterization}
In this section, we establish a combinatorial version of Afriat inequalities and deduce a potential-based characterization of rationalizability that eliminates the scaling multipliers $\blambda$. For any matrix $\BFM \in \mathbb{R}^{n \times n}$, let $\sign(\BFM) \in \{-1, 0, 1\}^{n \times n}$ denote its entrywise sign matrix, where $\sign(x) = 1$ if $x > 0$, $-1$ if $x < 0$, and $0$ if $x = 0$.
\begin{theorem}\label{thm:afr_ext}
Consider a consumption dataset $D$ with $n$ observations and the revealed preference matrix $\A \in \R^{n \times n}$ defined as in~\eqref{eq:revealedmatrix}. Let $\B \in \R^{n\times n}$ be any matrix with the same sign pattern as $\A$, i.e. $\sign(\B) = \sign(\A)$. The following statements are equivalent:
\begin{enumerate}[leftmargin=1cm]
    \item $\HM(D) = 1$;\label{cond1}
    \item $D$ satisfies GARP;\label{cond2}
    \item The following Afriat-type inequalities system $\AI(\B)$ is feasible:
    \begin{equation*}\label{sys:afr_ext}
 \bphi \in \mathbb{R}^n \quad \blambda \geq \mathbf{1} \quad \text{ and }  \phi_i - \phi_j + \lambda_i b_{ij} \ge 0 \quad \for i,j \in [n], 
    \end{equation*}
    \label{cond3}
        where $\mathbf{1}$ is the all-ones vector and $b_{ij}$ is the $(i,j)$ entry of $\B$. 
\end{enumerate} 
\end{theorem}

Theorem~\ref{thm:afr_ext} shows that rationalizability depends solely on the sign pattern of the revealed preference matrix: any matrix $\B$ satisfying $\sign(\B) = \sign(\A)$ yields an equivalent feasibility test $\AI(\B)$. This sign invariance allows us to address Question~\ref{eq:question} affirmatively. Specifically, selecting a standardized matrix $\tilde{\A} \in \mathbb{R}^{n \times n}$  eliminates $\blambda$ altogether, reducing consistency testing to a system of difference constraints in the potential vector $\bphi$.

% By carefully choosing a specific, bounded matrix $\tilde{\A} \in \R^{n \times n}$ such that $\sign(\tilde{\A}) = \sign(\A)$, we can completely eliminate $\blambda$ from the feasibility check. Our observation is formally stated in the following corollary, as a consequence of Theorem~\ref{thm:afr_ext}. We will provide a rigorous proof for both Theorem~\ref{thm:afr_ext} and the following corollary in this section, along with an intuitive geometric interpretation in Section~\ref{subsec:potential}.

\begin{corollary}[Potential-Based Characterization]\label{cor:potential}
Consider a consumption dataset $D$ and the revealed preference matrix $\A$. Define the standardized matrix $\tilde{\A} = (\tilde{a}_{ij}) \in \R^{n \times n}$ by
\[
\tilde{a}_{ij} := \frac{n}{n+1} \text{ if } a_{ij} > 0 \qquad \tilde{a}_{ij} := -\frac{1}{n+1} \text{ if } a_{ij} <0 \qquad \text{ and } \tilde{a}_{ij} =0 \text{ otherwise}. 
\]
% \begin{equation*}
%     \tilde{a}_{ij} :=
% \begin{cases}
%     \frac{n}{n+1} & \text{if } a_{ij} > 0, \\
%     -\frac{1}{n+1} & \text{if } a_{ij} < 0, \\
%     0 & \text{otherwise.}
% \end{cases}
% \end{equation*}
Then, $\HM(D) = 1$ if and only if there exists a vertex potential vector $\bphi \in \mathbb{R}^n$ satisfying
\begin{equation*}
\phi_j - \phi_i \le \tilde{a}_{ij} \quad \for  i, j \in [n].
\end{equation*}
\end{corollary}

\begin{figure}[htbp]
    \centering
    \begin{tikzpicture}[x=1cm,y=1cm,>=stealth,every node/.style={font=\small},line cap=round,line join=round]
        % Hatched negative orthant R^n_- in all three panels
        \foreach \panelshift in {0,5.50,11.00} {
            \begin{scope}[xshift=\panelshift cm]
                \clip (-1.35,-1.45) rectangle (0,0);
                \foreach \hatchpos in {-2.95,-2.70,...,0.05} {
                    \draw[black!45,line width=0.35pt] (\hatchpos,-1.45) -- ++(1.60,1.60);
                }
            \end{scope}
        }

        % Panel (b): Overlapping convex hulls conv(C_A) and conv(C_B)
        \begin{scope}[xshift=5.50cm]
            \filldraw[fill=black,fill opacity=0.11,draw=black!80,line width=0.65pt]
                (0,0) -- (1.55,0.80) -- (0.55,1.95) -- (-0.72,1.55) -- cycle;
            \filldraw[fill=black,fill opacity=0.11,draw=black!80,line width=0.65pt]
                (0,0) -- (2.50,0.42) -- (1.55,1.72) -- (-1.12,0.82) -- cycle;
        \end{scope}

        % Panel (c): Convex hull of C_B and separating hyperplane
        \begin{scope}[xshift=11.00cm]
            \filldraw[fill=black!12,draw=black,line width=0.65pt]
                (0,0) -- (2.50,0.42) -- (1.55,1.72) -- (-1.12,0.82) -- cycle;
            \draw[dashed,semithick] (-1.35,0.6075) -- (2.45,-1.1025);
            \draw[->,thick] (0,0) -- (0.45,1.00) node[right=1pt] {$\blambda$};
            \node[anchor=north east,font=\footnotesize] at (2.82,-1.14) {$\langle\blambda,\x\rangle=0$};
        \end{scope}

        % Axes and R^n_- label in all three panels
        \foreach \panelshift in {0,5.50,11.00} {
            \begin{scope}[xshift=\panelshift cm]
                \draw[->,semithick] (-1.45,0) -- (2.88,0);
                \draw[->,semithick] (0,-1.55) -- (0,2.75);
                \node[fill=white,inner sep=1.5pt] at (-0.68,-0.92) {$\R^n_-$};
            \end{scope}
        }

        % Panel (a): 4 points of C_A
        \begin{scope}[xshift=0cm]
            \foreach \pointx/\pointy in {-0.72/1.55,0.55/1.95,1.55/0.80} {
                \fill (\pointx,\pointy) circle (1.6pt);
            }
            \fill[white] (0,0) circle (3pt);
            \fill (0,0) circle (1.6pt);
            \node[anchor=north east,font=\footnotesize] at (2.75,2.48) {$\mathcal{C}_{\A}$};
            \node at (0.72,3.05) {(a) Cyclic embedding};
        \end{scope}

        % Panel (b): Sign invariance (C_A -> C_B and overlapping convex hulls)
        \begin{scope}[xshift=5.50cm]
            \draw[->,semithick,black!70,shorten <=2.8pt,shorten >=3.2pt] (-0.72,1.55) -- (-1.12,0.82);
            \draw[->,semithick,black!70,shorten <=2.8pt,shorten >=3.2pt] (0.55,1.95) -- (1.55,1.72);
            \draw[->,semithick,black!70,shorten <=2.8pt,shorten >=3.2pt] (1.55,0.80) -- (2.50,0.42);
            \foreach \pointx/\pointy in {-0.72/1.55,0.55/1.95,1.55/0.80} {
                \fill (\pointx,\pointy) circle (1.6pt);
            }
            \foreach \pointx/\pointy in {-1.12/0.82,1.55/1.72,2.50/0.42} {
                \filldraw[fill=white,draw=black,semithick] (\pointx,\pointy) circle (1.85pt);
            }
            \fill[white] (0,0) circle (3pt);
            \fill (0,0) circle (1.6pt);

            % Leader lines and labels for conv(C_A) and conv(C_B)
            \node[anchor=east,font=\footnotesize,inner sep=1.5pt] (lblA) at (-0.14,2.38) {$\conv(\mathcal{C}_{\A})$};
            \draw[thin,black!75] (lblA.south) -- (-0.36,1.66);
            \node[font=\footnotesize,inner sep=1.5pt] (lblB) at (1.90,-0.85) {$\conv(\mathcal{C}_{\B})$};
            \draw[thin,black!75] (lblB.north) -- (1.80,0.30);

            \node at (0.72,3.05) {(b) Sign invariance};
        \end{scope}

        % Panel (c): 4 points of C_B on conv(C_B)
        \begin{scope}[xshift=11.00cm]
            \foreach \pointx/\pointy in {-1.12/0.82,1.55/1.72,2.50/0.42} {
                \filldraw[fill=white,draw=black,semithick] (\pointx,\pointy) circle (1.85pt);
            }
            \fill[white] (0,0) circle (3pt);
            \filldraw[fill=white,draw=black,semithick] (0,0) circle (1.85pt);
            \node[anchor=north east,font=\footnotesize] at (2.75,2.48) {$\conv(\mathcal{C}_{\B})$};
            \node at (0.72,3.05) {(c) Convex hull separation};
        \end{scope}

    \end{tikzpicture}
    \caption{Proof architecture of Theorem~\ref{thm:afr_ext}. Filled and hollow circles denote $\mathcal{C}_{\A}$ and $\mathcal{C}_{\B}$, respectively, for matrices $\B$ satisfying $\sign(\B) = \sign(\A)$}
    \label{fig:thm2-proof}
\end{figure}
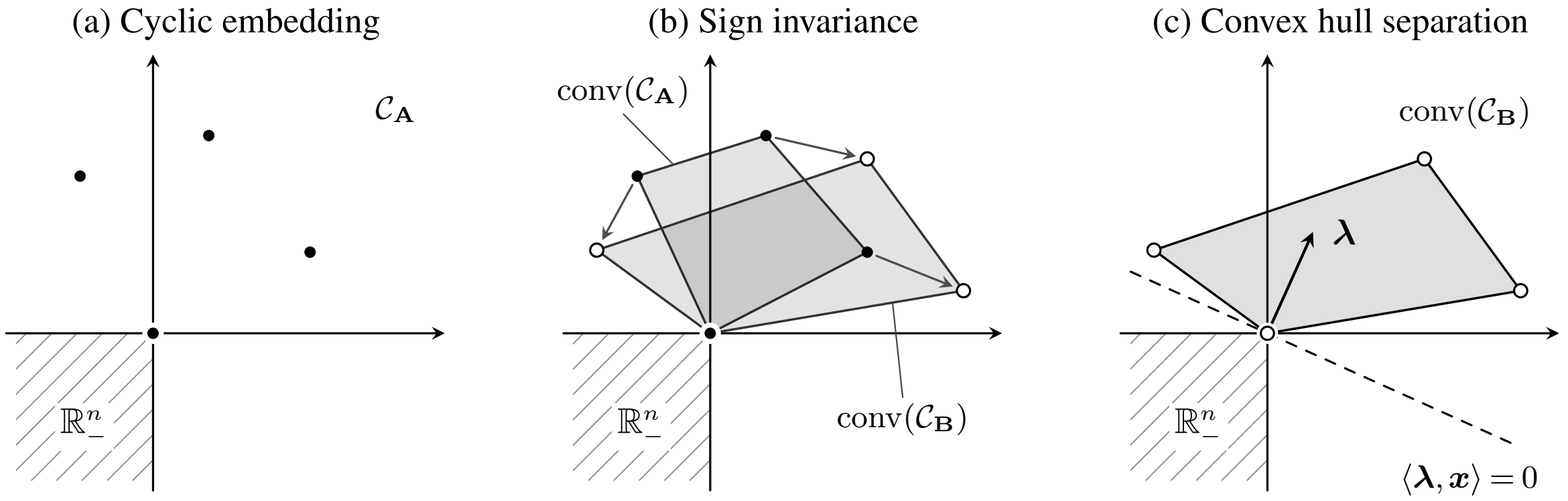

The remainder of this section is devoted to proving Theorem~\ref{thm:afr_ext} and Corollary~\ref{cor:potential}. For Theorem~\ref{thm:afr_ext}, the main challenge is showing that Condition~\ref{cond2} implies Condition~\ref{cond3}. We address this challenge in three steps, corresponding from left to right across the panels of Figure~\ref{fig:thm2-proof}. First, in Section~\ref{sec:embedding} (Panel~(a)), we construct the cyclic embedding set $\mathcal{C}_{\A}$ to encode cyclic comparisons in the choice data and show that $D$ satisfies GARP if and only if $\mathcal{C}_{\A} \cap \R^n_- = \emptyset$. Second, in Section~\ref{sec:orthant} (Panel~(b)), we prove that both $\conv(\mathcal{C}_{\A})$ and $\conv(\mathcal{C}_{\B})$ remain disjoint from $\R^n_-$, establishing sign invariance of this disjointness condition via the \textit{orthant hull}. Third, in Section~\ref{sec:separation} (Panel~(c)), we separate $\conv(\mathcal{C}_{\B})$ from $\R^n_-$ by a hyperplane to obtain a normal vector $\blambda \ge \mathbf{1}$, which yields a feasible solution to $\AI(\B)$. Finally, we complete the proof of Theorem~\ref{thm:afr_ext} in Section~\ref{sec:proof-thm} and prove Corollary~\ref{cor:potential} in Section~\ref{subsec:potential}.

\subsection{A disjointness Test via Cyclic Embedding}\label{sec:embedding}

Let $D = \{(\p^i, \q^i)\}_{i\in [n]}$ be a consumption dataset with revealed preference matrix $\A = (a_{ij}) \in \R^{n \times n}$.  By definition, $D$ violates GARP if and only if some cycle $\boldsymbol{c}$ satisfies $a_{c_\ell, c_{\ell+1}} \le 0$ for all $\ell \in [L]$ with at least one strict inequality, where a \emph{cycle} $\boldsymbol{c} = (c_1, \ldots, c_L, c_{L+1})$ of length $L \ge 1$ consists of distinct indices $c_1, \ldots, c_L \in [n]$ with $c_{L+1} = c_1$. We translate this GARP condition into a geometric disjointness test in $\R^n$, with each coordinate indexing an observation. In particular, we embed cycle weights to vectors as follows, which has been adopted in~\cite{beggs2021afriat}. 

\begin{definition}
    For a cycle $\boldsymbol{c} = (c_1, \ldots, c_L, c_{L+1})$ and matrix $\B \in \R^{n \times n}$, the cyclic embedding vector $\chi(\B, \boldsymbol{c}) \in \R^n$ has coordinates
    \begin{equation*}
        \chi_i(\B, \boldsymbol{c})  =  b_{i, c_{\ell+1}} \quad \text{if } i = c_\ell \text{ for some } \ell \in [L] \qquad  \text{ and } \qquad \chi_i(\B, \boldsymbol{c}) = 0 \text{ otherwise.}
    \end{equation*}
    The \emph{cyclic embedding set} is $\mathcal{C}_{\B} := \bigl\{ \chi(\B, \boldsymbol{c}) \bigm| \boldsymbol{c} \text{ is a cycle in } [n] \bigr\}$.
\end{definition}

% For $\B = \A$, each coordinate $\chi_i(\A, c)$ records $a_{c_\ell, c_{\ell+1}}$ when $i = c_\ell$ lies on $c$, and $0$ otherwise. Because the indices $c_1, \ldots, c_L$ are distinct, $\chi(\A, c)$ is well-defined and captures all pairwise comparisons along $c$. By definition, $D$ violates GARP if and only if some cycle $c$ satisfies $a_{c_\ell, c_{\ell+1}} \le 0$ for all $\ell \in [L]$ with at least one strict inequality—equivalently, $\chi(\A, c) \in \R^n_- := \R^n_{\le 0} \setminus \{\mathbf{0}\}$. Consequently, $D$ satisfies GARP if and only if
% \[
% \mathcal{C}_{\A} \cap \R^n_- = \emptyset,
% \]
% and since $\sign(\B) = \sign(\A)$, this holds if and only if $\mathcal{C}_{\B} \cap \R^n_- = \emptyset$.

By isolating the entries in the revealed preference matrix $\A$ along the active edges of a cycle and compressing the remaining coordinates to zero, $\mathcal{C}_{\A}$ maps the preference relations directly into the geometry of $\R^n$. Under this construction, a violation of rationalizability (a GARP violation) manifests when the embedding $\mathcal{C}_{\A}$ contains a nonzero vector in the negative orthant. This yields a geometric representation of GARP. Let $\mathbb{R}^n_- :=  \mathbb{R}^n_{\le 0} \setminus \{0\}$. Then $D$ satisfies GARP if and only if
\[
\mathcal{C}_{\A} \cap \R^n_{-} = \emptyset. 
\]
Moreover, this condition is sign invariant; that is, for any matrix $\B \in \R^{n \times n}$ satisfying $\sign(\B) = \sign(\A)$, GARP holds if and only if $\mathcal{C}_{\B} \cap \R^n_- = \emptyset$.

\subsection{Disjointness Preservation under the Orthant Hull}\label{sec:orthant}
In Section~\ref{sec:embedding}, we reduce GARP verification to the disjointness condition $\mathcal{C}_{\A} \cap \R^n_{-} = \emptyset$. This disjointness is preserved when taking the convex hull of $\mathcal{C}_{\A}$. This is an intriguing property of the cyclic embedding, as such preservation does not hold for general sets. More interestingly, this property is sign invariant: for any matrix $\B \in \R^{n \times n}$ satisfying $\sign(\B)= \sign(\A)$, $\conv(\mathcal{C}_{\B}) \cap \R^n_- = \emptyset$. In this subsection, we unify these phenomena through the \emph{orthant hull}. Informally, the orthant hull of $S \subseteq \R^n$ forms a convex relaxation while preserving the sign patterns present in $S$.

To define this operation formally, we extend the entrywise sign map to subsets $S \subseteq \R^n$ by $\sign(S) := \bigl\{\sign(\x) \in \{0,+,-\}^n \bigm| \x \in S\bigr\}$, where $\sign(\x) := \bigl(\sign(x_1), \ldots, \sign(x_n)\bigr)$. Conversely, for any subset of sign vectors $\mathcal{S} \subseteq \{0,+,-\}^n$, its preimage under $\sign$ is $\sign^{-1}(\mathcal{S}) := \bigl\{\x \in \R^{n} \bigm| \sign(\x) \in \mathcal{S} \bigr\}$. Using these maps, we define the orthant hull as follows.
\begin{definition}
    For a set $S \subseteq \R^n$, its \emph{orthant hull}, denoted by $\orth(S)$, is 
    \[
    \orth(S) := \conv \bigl(\sign^{-1}\bigl(\sign(S)\bigr)\bigr).
    \]
In words, the orthant hull of $S$ is the convex hull of the union of orthants generated by points in  $S$. 
\end{definition}

\begin{proposition}\label{prop:nonempty}
  Consider a consumption dataset $D$ and the revealed preference matrix $\A$. Let $\B$ be a matrix such that $\sign(\B) = \sign(\A)$. Then, $D$ satisfies GARP if and only if each of the following equivalent statements holds:
    \begin{enumerate}[leftmargin=1cm]
    \item  $\mathcal{C}_\A \cap \mathbb{R}_-^n  = \emptyset$, \label{prop:nonempty-1}
        \item  $\conv(\mathcal{C}_\A) \cap \mathbb{R}_-^n   = \emptyset$,\label{prop:nonempty-2}
        \item $\orth(\mathcal{C}_\A)  \cap  \mathbb{R}_-^n  = \emptyset$,\label{prop:nonempty-3}
        \item       $\orth(\mathcal{C}_\B)  \cap  \mathbb{R}_-^n  = \emptyset$.\label{prop:nonempty-4}
    \end{enumerate}
\end{proposition}

The proof exploits a combinatorial representation of the orthant hull. Namely, for any set $S \subseteq \R^n$, the sign patterns realized in $\orth(S)$ coincide with the span of $\Sign(S)$ in the sign space $\{0,+,-\}^n$. To define this span, we introduce the hypersum of sign vectors, which was introduced in matroid theory by \citet{baker2019matroids} to generalize oriented matroids; we follow the presentation of \citet{Anderson_2025}. Specifically, the hypersum $\boxplus$ on $\{0,+,-\}$ is defined by

\begin{equation*}
    + \boxplus + = \{ + \}
\end{equation*}
\begin{equation*}
    - \boxplus - = \{ - \}
\end{equation*}
\begin{equation*}
    X \boxplus 0 = 0 \boxplus X = \{X\} \quad \for  X \in \{0,+,-\}
\end{equation*}
\begin{equation*}
    + \boxplus - = - \boxplus + = \{0,+,-\}.
\end{equation*}
For any two sign vectors $\X, \Y \in \{0,+,-\}^n$, their coordinatewise hypersum is the set
\[
\X \boxplus \Y := \bigl\{ \Z \in \{0,+,-\}^n \bigm| Z_i \in X_i \boxplus Y_i \ \for i \in [n] \bigr\}.
\]
Because $\boxplus$ is commutative and associative, for any nonempty finite subset of sign vectors $\mathcal{S} = \{\X^1, \ldots, \X^k\} \subseteq \{0,+,-\}^n$, we write $\boxplus_{\X \in \mathcal{S}} \X$ for the iterated hypersum $\X^1 \boxplus \cdots \boxplus \X^k$, with $\boxplus_{\X \in \{\X^1\}} \X := \{\X^1\}$ for a singleton.

\begin{definition}
    For $\mathcal{S}\subseteq\{0,+,-\}^n$, the span of $\mathcal{S}$, denoted by $\Span(\mathcal{S})$, is the set of all sign vectors that belong to the hypersum of some nonempty subset of $\mathcal{S}$. Equivalently,
\begin{equation*}
    \Span(\mathcal{S}) = \bigcup_{\varnothing\ne\mathcal{T}\subseteq\mathcal{S}} \left(\boxplus_{\X\in\mathcal{T}}\X\right).
\end{equation*}
\end{definition}
Equivalently, $\Span(\mathcal{S})$ is the minimal subset of $\{0,+,-\}^n$ containing $\mathcal{S}$ that is closed under hypersum addition.

With the sign span defined, we can now formalize the exact correspondence between continuous convexification and discrete hypersum closure. Specifically, taking the convex hull of a union of sign chambers in $\R^n$ commutes with the preimage map under the sign operator, reducing to taking the hypersum span in $\{0,+,-\}^n$. Our key structural lemma makes this equivalence precise and gives a direct combinatorial formula for the orthant hull.

\begin{lemma}\label{lemma:orth-span}
    For any set of sign vectors $\mathcal{S} \subseteq  \{0,+,-\}^n$, the convex hull of the inverse image of $\mathcal{S}$ under $\sign$ equals the inverse image of its hypersum span under $\sign$: 
    \begin{equation*}
        \conv\bigl(\sign^{-1}(\mathcal{S})\bigr) = \sign^{-1}\bigl(\Span(\mathcal{S})\bigr).
    \end{equation*}
    In particular, for a set $S \subseteq \R^n$, $\Orth (S) = \Sign^{-1}\bigl(\Span(\Sign(S)) \bigr)$. 
\end{lemma}
\proof{Proof of Lemma~\ref{lemma:orth-span}.}
We first show that $\conv\bigl(\sign^{-1}(\mathcal{S})\bigr) \subseteq \sign^{-1}\bigl(\Span(\mathcal{S})\bigr)$. Since $\mathcal{S} \subseteq \Span(\mathcal{S})$, we have $\sign^{-1}(\mathcal{S}) \subseteq \sign^{-1}(\Span(\mathcal{S}))$. It therefore suffices to show that $\sign^{-1}(\Span(\mathcal{S}))$ is convex. Take arbitrary $\x, \boldsymbol{y} \in \sign^{-1}(\Span(\mathcal{S}))$ and let $\boldsymbol{z} := \lambda \x + (1 - \lambda)\boldsymbol{y}$ for some $\lambda \in (0, 1)$. Set $\X := \sign(\x)$ and $\Y := \sign(\boldsymbol{y})$. We claim that $\sign(z_i) \in X_i \boxplus Y_i$ for every $i \in [n]$. Indeed, if $X_i = Y_i$, then $z_i$ has that same sign (and is zero if both are zero); if exactly one of $X_i, Y_i$ is zero, then $z_i$ takes the sign of the nonzero term because both $\lambda$ and $1 - \lambda$ are strictly positive; if $X_i$ and $Y_i$ have opposite signs, then $X_i \boxplus Y_i = \{0, +, -\}$, which contains every possible sign of $z_i$. Hence $\sign(\boldsymbol{z}) \in \X \boxplus \Y$. Because $\X, \Y \in \Span(\mathcal{S})$ and $\Span(\mathcal{S})$ is closed under hypersums, we have $\X \boxplus \Y \subseteq \Span(\mathcal{S})$, and thus $\sign(\boldsymbol{z}) \in \Span(\mathcal{S})$. It follows that $\boldsymbol{z} \in \sign^{-1}(\Span(\mathcal{S}))$, proving convexity of $\sign^{-1}(\Span(\mathcal{S}))$ and hence $\conv\bigl(\sign^{-1}(\mathcal{S})\bigr) \subseteq \sign^{-1}\bigl(\Span(\mathcal{S})\bigr)$.

For the reverse inclusion, take $\boldsymbol{z}\in\sign^{-1}(\Span(\mathcal{S}))$. By the definition of the span, there exist $\X^1,\ldots,\X^k\in\mathcal{S}$ such that $\sign(\boldsymbol{z})\in\X^1\boxplus\cdots\boxplus\X^k$. We will construct vectors $\x^1,\ldots,\x^k$ with $\sign(\x^j)=\X^j \in \mathcal{S}$ for each $j$ and $\boldsymbol{z}=\frac{1}{k}\sum_{j=1}^k\x^j$, which will lead to $\z \in \conv(\sign^{-1}(\mathcal{S}))$ since $\x^j \in \sign^{-1}(\mathcal{S})$ and $\z$ is a convex combination of vectors in $\sign^{-1}(\mathcal{S})$.

Fix a coordinate $i$. We construct $x_i^1,\ldots,x_i^k$ such that  $\Sign(x_i^j) = X^j_i$ for $j \in [k]$ and $\sum_{j=1}^kx_i^j=kz_i$. Let $r$ and $s$ denote the numbers of positive and negative signs among $X_i^1,\ldots,X_i^k$, respectively. For $j \in [k]$ such that $X_i^j=0$, we set $x_i^j=0$. For $j \in [k]$ such that $X_i^j\neq 0$, we discuss three cases: 

\textbf{Case 1:} If $z_i>0$, then $+\in X_i^1\boxplus\cdots\boxplus X_i^k$. If $r=0$, every $X_i^j$ would be either $0$ or $-$, which is impossible since their hypersum would then be $\{0\}$ or $\{-\}$. Thus $r>0$. For $j \in [k]$ such that $X_i^j\neq 0$, let
\begin{equation*}
    x_i^j =\frac{kz_i+s}{r} \quad \text{ if } X_i^j=+ \quad \text{ and } \quad x_i^j  =-1 \quad \text{if } X_i^j=-. 
\end{equation*}
Therefore, $\sum_{j=1}^k x_i^j = r\frac{kz_i+s}{r}-s = kz_i$.

\textbf{Case 2:} If $z_i<0$, then $-\in X_i^1\boxplus\cdots\boxplus X_i^k$. If $s=0$, the hypersum would be $\{0\}$ if all $X_i^j$ were zero and would be $\{+\}$ otherwise. Thus $s>0$.  For $j \in [k]$ such that $X_i^j\neq 0$, let
\begin{equation*}
    x_i^j = 1  \quad \text{if } X_i^j=+ \quad \text{ and } \quad x_i^j = \frac{kz_i-r}{s}  \quad \text{if } X_i^j=-. 
\end{equation*}
Therefore, $\sum_{j=1}^k x_i^j = r+s\frac{kz_i-r}{s} = kz_i$. 

\textbf{Case 3:} If $z_i=0$, then $0\in X_i^1\boxplus\cdots\boxplus X_i^k$. If $r=s=0$, all $x_i^j$ have already been set to zero. Otherwise, both signs must occur, since the hypersum would be $\{+\}$ if only positive signs and zeros occurred, and would be $\{-\}$ if only negative signs and zeros occurred. In this case, set
\begin{equation*}
    x_i^j =  \frac{1}{r} \quad \text{if }  X_i^j=+ \quad \text{ and } \quad  x_i^j =  - \frac{1}{s} \quad \text{if } X_i^j=-.
\end{equation*}
Therefore, $\sum_{j=1}^k x_i^j = r\frac{1}{r}-s\frac{1}{s} = 0 = kz_i$.

Carrying out this construction in every coordinate gives $\sign(\x^j)=\X^j\in\mathcal{S}$ for every $j \in [k]$ and $\sum_{j=1}^k\x^j=k\boldsymbol{z}$. Thus, $\boldsymbol{z}$ is the average of $k$ vectors in $\sign^{-1}(\mathcal{S})$ and belongs to $\conv(\sign^{-1}(\mathcal{S}))$, showing that $\boldsymbol{z}\in\sign^{-1}(\Span(\mathcal{S})) \subseteq \conv(\sign^{-1}(\mathcal{S}))$. 

Last,  we can conclude that $\orth(S) = \conv(\sign^{-1}(\sign(S))) = \sign^{-1}(\Span(\sign(S)))$, where the first equality holds by definition, and the second equality has been shown above. \Halmos
\endproof

With these preparations, we can now prove Proposition~\ref{prop:nonempty}.

\proof{Proof of Proposition~\ref{prop:nonempty}}
By Section~\ref{sec:embedding}, $D$ satisfies GARP if and only if Condition~\ref{prop:nonempty-1} holds. The inclusions $\mathcal{C}_{\A}\subseteq\conv(\mathcal{C}_{\A})\subseteq\orth(\mathcal{C}_{\A})$ give the implications Condition $\ref{prop:nonempty-3}\Rightarrow \text{Condition } \ref{prop:nonempty-2}\Rightarrow  \text{Condition } \ref{prop:nonempty-1}$. Moreover, $\Sign(\A) = \Sign(\B)$ implies that  $\chi(\A,c)$ and $\chi(\B,c)$ have the same sign pattern for every cycle $c$, showing $\sign(\mathcal{C}_{\A})=\sign(\mathcal{C}_{\B})$. According to the definition of the orthant hull, $\orth(\mathcal{C}_{\A})=\orth(\mathcal{C}_{\B})$, proving the equivalence of Conditions~\ref{prop:nonempty-3} and~\ref{prop:nonempty-4}.

It remains to prove that Condition~\ref{prop:nonempty-1} implies Condition~\ref{prop:nonempty-3}. Suppose, to the contrary, that Condition~\ref{prop:nonempty-1} holds but there exists $\z\in\orth(\mathcal{C}_{\A})\cap\R^n_-$.  We will use $\z$ to construct a cycle that violates Condition~\ref{prop:nonempty-1}.

Let $\Z:=\sign(\z)\in\Span(\sign(\mathcal{C}_{\A}))$.   Since $\z\in\orth(\mathcal{C}_{\A})$ which, by Lemma~\ref{lemma:orth-span}, equals to  $\sign^{-1}(\Span(\sign(\mathcal{C}_{\A})))$, the definition of $\sign^{-1}$ gives $\Z \in\Span(\sign(\mathcal{C}_{\A}))$. By the definition of the span, there exist $\X^1,\ldots,\X^k\in\sign(\mathcal{C}_{\A})$ such that $\Z\in\X^1\boxplus\cdots\boxplus\X^k$. For each $j\in[k]$, membership in $\sign(\mathcal{C}_{\A})$ means that $\X^j=\sign(\chi(\A, \boldsymbol{c}^j))$ for some cycle $\boldsymbol{c}^j$. From these cycles, we will a cycle violating Condition~\ref{prop:nonempty-1}. 

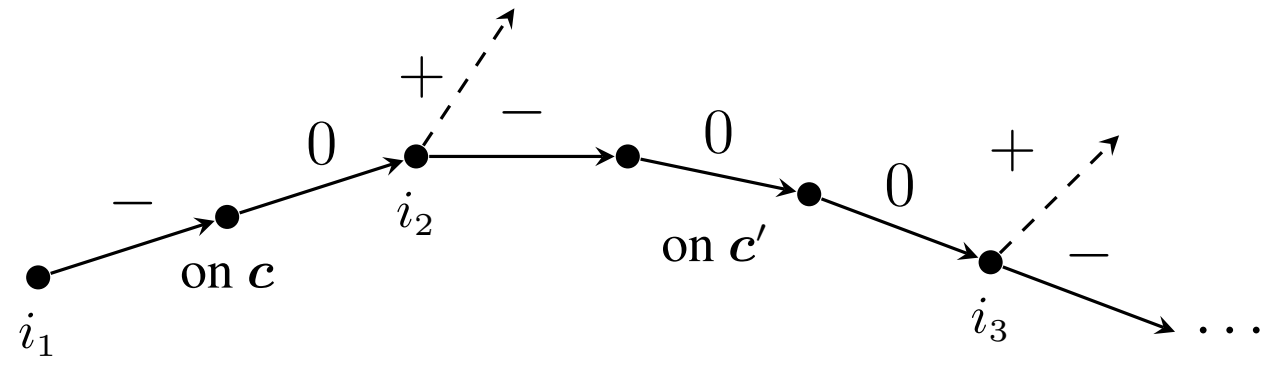
\begin{figure}[htbp]
    \centering
    \begin{tikzpicture}[>=stealth, vertex/.style={circle,fill,inner sep=1.6pt}, taken/.style={->,semithick}, skipped/.style={->,semithick,dashed}, every label/.style={font=\small}]
        \node[vertex,label=below:$i_1$] (a) at (0,0) {};
        \node[vertex] (b) at (1.25,0.40) {};
        \node[vertex,label=below:$i_2$] (d) at (2.50,0.80) {};
        \node[vertex] (e) at (3.90,0.80) {};
        \node[vertex] (f) at (5.10,0.55) {};
        \node[vertex,label=below:$i_3$] (g) at (6.30,0.10) {};

        \draw[taken] (a) -- node[above] {$-$} (b);
        \draw[taken] (b) -- node[above] {$0$} (d);
        \draw[skipped] (d) -- node[left] {$+$} (3.15,1.78);
        \draw[taken] (d) -- node[above] {$-$} (e);
        \draw[taken] (e) -- node[above] {$0$} (f);
        \draw[taken] (f) -- node[above] {$0$} (g);
        \draw[skipped] (g) -- node[above left] {$+$} (7.15,0.94);
        \draw[taken] (g) -- node[above] {$-$} (7.52,-0.36);
        \node[right] at (7.50,-0.36) {$\cdots$};

        \node[below] at (1.25,0.24) {\small on $\boldsymbol{c}$};
        \node[below] at (4.48,0.50) {\small on $\boldsymbol{c}'$};
    \end{tikzpicture}
    \caption{The walk follows the solid transitions. At $i_2$, it stops before the dashed positive transition on $\boldsymbol{c}$ and switches to the negative transition on $\boldsymbol{c}'$; the same switch occurs at $i_3$.}
    \label{fig:cycle-switch}
\end{figure}
\vspace{-0.3cm}

Because $\boldsymbol{z} \in \R^n_-$, every coordinate of $\Z$ lies in $\{0, -\}$, and $Z_{i_1} = -$ for some index $i_1 \in [n]$. Coordinatewise hypersum addition gives $Z_{i_1} \in X^1_{i_1} \boxplus \cdots \boxplus X^k_{i_1}$. If none of $X^1_{i_1}, \ldots, X^k_{i_1}$ equals $-$, their hypersum would be $\{0\}$ when all entries are $0$ and $\{+\}$ otherwise, neither of which contains $-$. Thus $X^j_{i_1} = -$ for some $j \in [k]$. Letting $\boldsymbol{c} := \boldsymbol{c}^j$, we have $\chi_{i_1}(\A, \boldsymbol{c}) < 0$, so $\boldsymbol{c}$ contains $i_1$ and its outgoing transition from $i_1$ has negative weight.

Start at $i_1$ and step forward along $\boldsymbol{c}$. Because $\boldsymbol{c}$ has a negative transition and Condition~\ref{prop:nonempty-1} holds, $\boldsymbol{c}$ must also contain at least one positive transition; see Figure~\ref{fig:cycle-switch}. Let $i_2$ be the first index encountered along $\boldsymbol{c}$ (starting from $i_1$) whose outgoing transition in $\boldsymbol{c}$ is positive, and stop at $i_2$ before taking that positive transition. Every transition traversed so far has sign $0$ or $-$, and the first step from $i_1$ has sign $-$.

To continue past $i_2$, observe that the outgoing transition of $\boldsymbol{c}$ at $i_2$ is positive, so $X^j_{i_2} = +$. Meanwhile, $Z_{i_2} \in \{0, -\}$ and $Z_{i_2} \in X^1_{i_2} \boxplus \cdots \boxplus X^k_{i_2}$. Since $X^j_{i_2} = +$, if no $X^{j'}_{i_2}$ were equal to $-$, the hypersum of $X^1_{i_2}, \ldots, X^k_{i_2}$ would be $\{+\}$, contradicting $Z_{i_2} \in \{0, -\}$. Hence $X^{j'}_{i_2} = -$ for some $j' \ne j$. Setting $\boldsymbol{c}' := \boldsymbol{c}^{j'}$, the cycle $\boldsymbol{c}'$ contains $i_2$ with a negative outgoing transition. We take this negative transition on $\boldsymbol{c}'$ and follow $\boldsymbol{c}'$ until its first positive transition at some index $i_3$, stopping before taking that positive transition; again, such a positive transition exists by Condition~\ref{prop:nonempty-1}.  Applying the same hypersum argument at $i_3$ yields another cycle with a negative outgoing transition from $i_3$.

Repeating this procedure produces an infinite walk of transitions in $[n]$, each of sign $0$ or $-$. Because $[n]$ is finite, some index must be visited at least twice. Extracting the closed segment between the first repeated visit and its next occurrence (removing any sub-loops if needed) yields a cycle $\boldsymbol{c}^\ast$ of distinct indices in $[n]$ whose transitions all have sign $0$ or $-$. Moreover, $\boldsymbol{c}^\ast$ includes at least one index at which the walk switched cycles and took a strictly negative transition. Hence $\chi(\A, \boldsymbol{c}^\ast) \in \mathcal{C}_{\A} \cap  \R^n_-$, contradicting Condition~\ref{prop:nonempty-1}. \hfill \Halmos
\endproof

\subsection{From Disjointness to Supporting Hyperplanes and Afriat Inequalities}\label{sec:separation}
When the convex hull of a cyclic embedding set is disjoint from $\R^n_{-}$, we can separate it using a supporting hyperplane with a normal vector $\blambda \ge \mathbf{1}$. This separation corresponds to the nonnegativity of all cycle lengths, whose dual yields Afriat inequalities.

\begin{proposition}\label{prop:disjointtoAI}
Let $\B \in \R^{n\times n}$ a matrix such that $\conv(\mathcal{C}_{\B})\cap \R^n_{-} = \emptyset $. Then, $\AI(\B)$ is feasible. 
\end{proposition}
\proof{Proof.}
We first show that there exists $\blambda \ge \mathbf{1}$ such that $\langle \blambda, \x \rangle \ge 0$ for every $\x \in \conv(\mathcal{C}_{\B})$. Let $K := \cone(\mathcal{C}_{\B})$ denote the conic hull of $\mathcal{C}_{\B}$ (the set of all nonnegative linear combinations of vectors in $\mathcal{C}_{\B}$). We claim that $\conv(\mathcal{C}_{\B}) \cap \R^n_{-} = \emptyset$ implies $K \cap \R^n_- = \emptyset$. Suppose to the contrary that there exists $\x \in K \cap \R^n_-$. Since $\x \ne \mathbf{0}$ and $\x \in K$, there exists a finite subset $\mathcal{P} \subseteq \mathcal{C}_{\B}$ and positive weights $\{w_p\}_{p \in \mathcal{P}}$ such that $\x = \sum_{p \in \mathcal{P}} w_p \boldsymbol{p}$. Setting $w := \sum_{p \in \mathcal{P}} w_p > 0$, we have $\frac{\x}{w} \in \R^n_-$ and
\[
\frac{\x}{w} = \sum_{p \in \mathcal{P}} \frac{w_p}{w} \boldsymbol{p} \in \conv(\mathcal{C}_{\B}),
\]
contradicting $\conv(\mathcal{C}_{\B}) \cap \R^n_- = \emptyset$.
Next, define $\Delta := \bigl\{\boldsymbol{y} \in \R^n_{\le 0} \bigm| \sum_{i=1}^n y_i = -1\bigr\}$. Since $\Delta \subseteq \R^n_-$ and $K \cap \R^n_- = \emptyset$, we have $K \cap \Delta = \emptyset$. Because $\mathcal{C}_{\B}$ is finite, $K$ is a closed convex cone, while $\Delta$ is compact and convex. By strong hyperplane separation (see, e.g., Corollary~11.4.2 and Theorem~11.1 of~\cite{rockafellar1970convex}), there exists a vector $\tilde{\boldsymbol{\lambda}} \in \R^n$ such that
\[
\sup_{\boldsymbol{y} \in \Delta} \langle \tilde{\boldsymbol{\lambda}}, \boldsymbol{y} \rangle < \inf_{\x \in K} \langle \tilde{\boldsymbol{\lambda}}, \x \rangle.
\]
If $\langle \tilde{\boldsymbol{\lambda}}, \x \rangle < 0$ for some $\x \in K$, then $t\x \in K$ for every $t > 0$ and $\langle \tilde{\boldsymbol{\lambda}}, t\x \rangle \to -\infty$ as $t \to \infty$, contradicting $\inf_{\x \in K} \langle \tilde{\boldsymbol{\lambda}}, \x \rangle > -\infty$. Hence $\langle \tilde{\boldsymbol{\lambda}}, \x \rangle \ge 0$ for all $\x \in K$. Since $\mathbf{0} \in K$, we obtain $\inf_{\x \in K} \langle \tilde{\boldsymbol{\lambda}}, \x \rangle = 0$, and therefore $\langle \tilde{\boldsymbol{\lambda}}, \boldsymbol{y} \rangle < 0$ for all $\boldsymbol{y} \in \Delta$. For each standard basis vector $\boldsymbol{e}_i$, we have $-\boldsymbol{e}_i \in \Delta$, so $\langle \tilde{\boldsymbol{\lambda}}, -\boldsymbol{e}_i \rangle = -\tilde{\lambda}_i < 0$, which gives $\tilde{\lambda}_i > 0$ for all $i \in [n]$. Rescaling by $\alpha := \bigl(\min_{i \in [n]} \tilde{\lambda}_i\bigr)^{-1} > 0$ and setting $\blambda := \alpha \tilde{\boldsymbol{\lambda}}$ yields $\blambda \ge \mathbf{1}$ while preserving $\langle \blambda, \x \rangle \ge 0$ for all $\x \in K$.

We now establish the feasibility of $\AI(\B)$. For any cycle $\boldsymbol{c} = (c_1, \ldots, c_L, c_{L+1})$ in $[n]$, the definition of the cyclic embedding vector and $\mathcal{C}_{\B} \subseteq K$ give
\[
\sum_{\ell=1}^{L} \lambda_{c_\ell} b_{c_\ell, c_{\ell+1}} = \bigl\langle \blambda, \chi(\B, \boldsymbol{c}) \bigr\rangle \ge 0.
\]
Thus every cycle in $[n]$ has nonnegative total weight under $w_{ij} := \lambda_i b_{ij}$. By potential feasibility on cycles (see~\cite{gallai1958maximum}), there exists $\bphi \in \R^n$ such that $\phi_j - \phi_i \le \lambda_i b_{ij}$ for all $i, j \in [n]$. Since $\blambda \ge \mathbf{1}$, the pair $(\bphi, \blambda)$ is feasible for $\AI(\B)$.\hfill \Halmos
\endproof

\subsection{Proof of  Combinatorial Afriat Characterization}\label{sec:proof-thm}
\proof{Proof of Theorem~\ref{thm:afr_ext}.}

We first show that Condition~\ref{cond2} implies Condition~\ref{cond3}. By Proposition~\ref{prop:nonempty}, if $D$ satisfies GARP, then $\orth(\mathcal{C}_{\B}) \cap \R^n_- = \emptyset$, and hence $\conv(\mathcal{C}_{\B}) \cap \R^n_- = \emptyset$. By Proposition~\ref{prop:disjointtoAI}, $\AI(\B)$ is feasible.

We next prove that Condition~\ref{cond1} implies Condition~\ref{cond2}. Suppose $\HM(D) = 1$, so a strictly increasing utility function $u : \R_{\ge 0}^m \to \R$ rationalizes $D$. Let $(\p^{i_1}, \q^{i_1}), \ldots, (\p^{i_L}, \q^{i_L})$ be any sequence of observations with $\langle \p^{i_k}, \q^{i_k} \rangle \ge \langle \p^{i_k}, \q^{i_{k+1}} \rangle$ for all $k \in [L]$ (setting $i_{L+1} := i_1$). By utility maximization,
\[
u(\q^{i_1}) \ge u(\q^{i_2}) \ge \cdots \ge u(\q^{i_L}) \ge u(\q^{i_1}),
\]
forcing $u(\q^{i_1}) = \cdots = u(\q^{i_L})$. If $\langle \p^{i_k}, \q^{i_k} \rangle > \langle \p^{i_k}, \q^{i_{k+1}} \rangle$ for some $k \in [L]$, then since $\p^{i_k} \in \R_{>0}^m$, there exists $\epsilon > 0$ such that $\q^{i_{k+1}} + \epsilon\mathbf{1}$ is affordable at prices $\p^{i_k}$, and strict monotonicity gives $u(\q^{i_{k+1}} + \epsilon\mathbf{1}) > u(\q^{i_{k+1}}) = u(\q^{i_k})$, contradicting optimality of $\q^{i_k}$. Hence $\langle \p^{i_k}, \q^{i_k} \rangle = \langle \p^{i_k}, \q^{i_{k+1}} \rangle$ for all $k \in [L]$, so $D$ satisfies GARP.

Finally, we show that Condition~\ref{cond3} implies Condition~\ref{cond1}. Suppose that $(\bphi, \blambda)$ is feasible for $\AI(\B)$. Summing the inequalities $\phi_i - \phi_j + \lambda_i b_{ij} \ge 0$ around any cycle $\boldsymbol{c}$ in $[n]$ gives $\langle \blambda, \chi(\B, \boldsymbol{c}) \rangle \ge 0$. Because $\lambda_i > 0$ for every $i \in [n]$, we have $\langle \blambda, \x \rangle < 0$ for every $\x \in \R^n_-$, and thus $\mathcal{C}_{\B} \cap \R^n_- = \emptyset$. Moreover, $\sign(\B) = \sign(\A)$ implies $\sign(\chi(\A, \boldsymbol{c})) = \sign(\chi(\B, \boldsymbol{c}))$, so $\chi(\A, \boldsymbol{c}) \in \R^n_-$ if and only if $\chi(\B, \boldsymbol{c}) \in \R^n_-$. Consequently, $\mathcal{C}_{\A} \cap \R^n_- = \emptyset$, and Proposition~\ref{prop:nonempty} implies that $D$ satisfies GARP. Applying the implication from Condition~\ref{cond2} to Condition~\ref{cond3} with $\B = \A$, we obtain a feasible pair $(\bphi, \blambda)$ for $\AI(\A)$.  Define $u : \R^m_{\ge 0} \to \R$ by
\[
u(\q) := \min_{i \in [n]} \bigl\{ \phi_i + \lambda_i \langle \p^i, \q - \q^i \rangle \bigr\}.
\]
As a pointwise minimum of strictly increasing affine functions ($\lambda_i > 0$, $\p^i \in \R^m_{>0}$), $u$ is concave and strictly increasing. For every $j \in [n]$, feasibility in $\AI(\A)$ gives $\phi_i + \lambda_i \langle \p^i, \q^j - \q^i \rangle = \phi_i + \lambda_i a_{ij} \ge \phi_j$ for all $i \in [n]$ (with equality at $i = j$), so $u(\q^j) = \phi_j$. Finally, for any $i \in [n]$ and any $\q \in \R^m_{\ge 0}$ with $\langle \p^i, \q \rangle \le \langle \p^i, \q^i \rangle$, we have $\langle \p^i, \q - \q^i \rangle \le 0$ and $\lambda_i > 0$, whence
\[
u(\q) \le \phi_i + \lambda_i \langle \p^i, \q - \q^i \rangle \le \phi_i = u(\q^i).
\]
Thus $u$ rationalizes $D$, and $\HM(D) = 1$. \hfill \Halmos
\endproof

\subsection{Proof of the Potential-Based Characterization}\label{subsec:potential}
Proposition~\ref{prop:nonempty} gives us the freedom to reposition $\mathcal{C}_{\B}$: because $\orth(\mathcal{C}_{\B}) = \orth(\mathcal{C}_{\A})$ whenever $\sign(\B) = \sign(\A)$, varying the magnitudes of the entries preserves both the geometric characterization of GARP and the separation of $\conv(\mathcal{C}_{\B})$ from $\R^n_-$. In this subsection, we exploit this invariance to choose a representation for which $\blambda$ can be fixed to $\mathbf{1} = (1, \ldots, 1)$. Fixing $\blambda$ eliminates the multiplier variables from $\AI(\B)$, yielding the potential-based characterization in Corollary~\ref{cor:potential}.

\begin{figure}[H]
    \centering
    \begin{tikzpicture}[scale=0.8,x=1.18cm,y=1.18cm,>=stealth,line cap=round,line join=round,every node/.style={font=\small}]
        \foreach \panelshift in {0,6.1} {
            \begin{scope}[xshift=\panelshift cm]
                \fill[black!4] (0,0) rectangle (2.3,2.65);
                \begin{scope}
                    \clip (-2,-1.35) rectangle (0,0);
                    \foreach \hatchpos in {-3.4,-3.15,...,0.1} {
                        \draw[black!35,line width=0.3pt] (\hatchpos,-1.35) -- ++(1.35,1.35);
                    }
                \end{scope}
                \node[fill=white,inner sep=1pt] at (-1.25,-0.95) {$\R^n_-$};
                \node at (1.75,2.5) {$\R^n_{\ge0}$};
            \end{scope}
        }

        \filldraw[fill=black!13,draw=black,semithick] (0,0) -- (-1.8,0.6) -- (0.8,1.5) -- (1.2,0.8) -- cycle;
        \begin{scope}[xshift=6.1cm]
            \filldraw[fill=black!13,draw=black,semithick] (0,0) -- (-0.36,1.65) -- (0.8,1.5) -- (1.2,0.8) -- cycle;
        \end{scope}

        \foreach \panelshift in {0,6.1} {
            \begin{scope}[xshift=\panelshift cm]
                \draw[->,semithick] (-2.05,0) -- (2.45,0);
                \draw[->,semithick] (0,-1.4) -- (0,2.8);
            \end{scope}
        }

        \draw[dashed,semithick] (-2,0.5) -- (1.7,-0.425);
        \draw[->,thick] (0,0) -- (0.3,1.2) node[right] {$\blambda$};
        \foreach \pointx/\pointy in {-1.8/0.6,0.8/1.5,1.2/0.8,0/0} {
            \fill (\pointx,\pointy) circle (1.6pt);
        }
        \node at (0.15,3.15) {(a) An arbitrary position};
        \node[anchor=east] at (-0.2,2.25) {$\conv(\mathcal{C}_{\B})$};

        \begin{scope}[xshift=6.1cm]
            \draw[densely dotted,thick] (-1.65,1.65) -- (1.3,-1.3);
            \draw[->,thick] (0,0) -- (0.65,0.65) node[above left] {$\mathbf{1}$};
            \draw[->,semithick,shorten <=3pt,shorten >=3pt] (-1.8,0.6) -- (-0.36,1.65);
            \draw[fill=white,semithick] (-1.8,0.6) circle (1.8pt);
            \foreach \pointx/\pointy in {-0.36/1.65,0.8/1.5,1.2/0.8,0/0} {
                \fill (\pointx,\pointy) circle (1.6pt);
            }
            \node[anchor=north east] at (2.2,-1.35) {$\langle\mathbf{1},\x\rangle=0$};
            \node at (0.15,3.15) {(b) A fixed separating normal};
            \node[anchor=east] at (-0.2,2.25) {$\conv(\mathcal{C}_{\B'})$};
        \end{scope}
    \end{tikzpicture}
    \caption{The vertex within $\{ \x \mid \langle\mathbf{1},\x\rangle<0 \}$ moves upward and to the right: its positive coordinate increases and its negative coordinate moves toward zero. The hollow circle marks the original position; after the move, $\mathbf{1}$ separates $\conv(\mathcal{C}_{\B'})$ from $\R^n_-$.}
    \label{fig:fixed-normal-separation}
\end{figure}
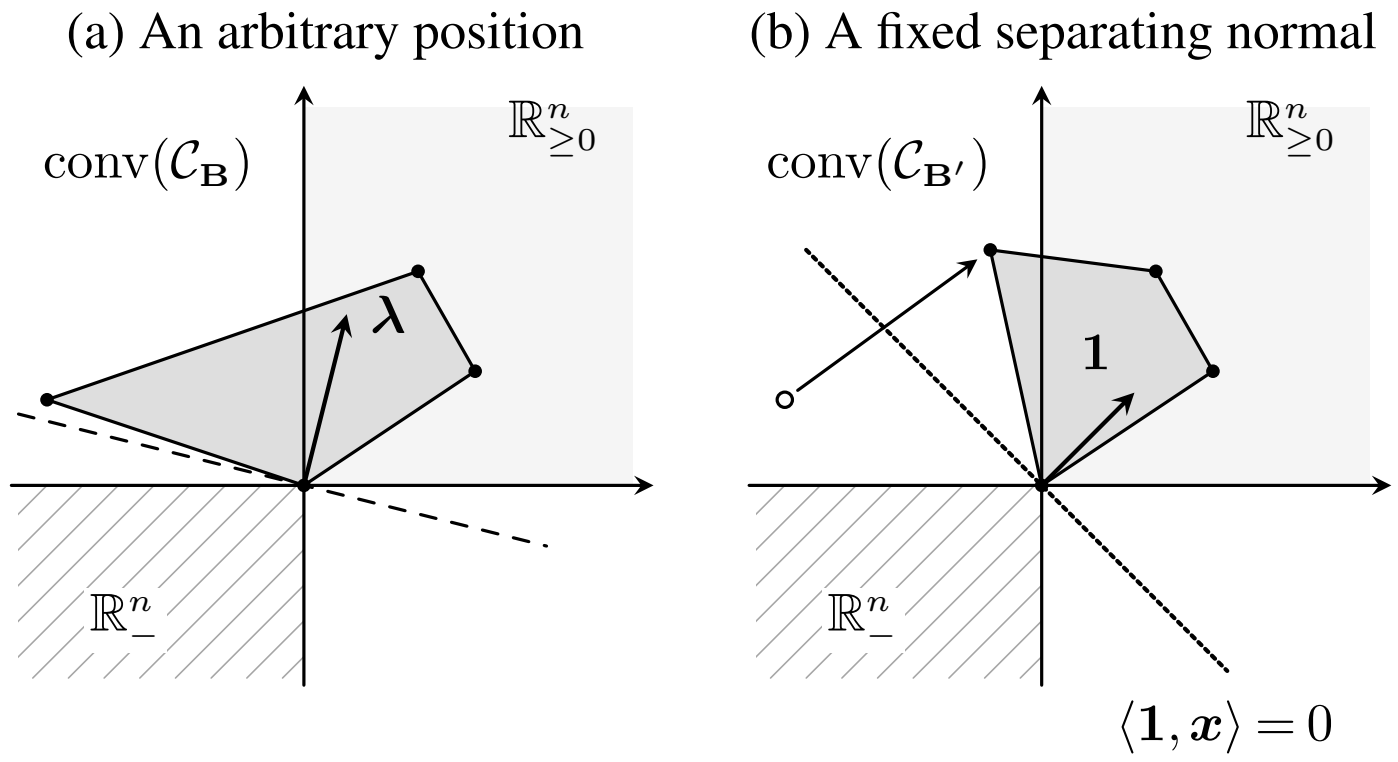

In the proof of Theorem~\ref{thm:afr_ext}, $\blambda$ acts as a positive normal vector separating $\conv(\mathcal{C}_{\B})$ from $\R^n_-$. To fix $\blambda = \mathbf{1}$, the set $\mathcal{C}_{\B}$ must be positioned so that the hyperplane $\langle \mathbf{1}, \x \rangle = 0$ already separates it from $\R^n_-$, meaning
\begin{equation*}
    \langle \mathbf{1}, \x \rangle \ge 0 \qquad \text{for every } \x \in \mathcal{C}_{\B}.
\end{equation*}
Geometrically, the freedom granted by Proposition~\ref{prop:nonempty} lets us rotate $\mathcal{C}_{\B}$ toward the nonnegative orthant $\R^n_{\ge 0}$ until the fixed normal $\mathbf{1}$ suffices, as illustrated in Figure~\ref{fig:fixed-normal-separation}. Algebraically, this corresponds to making the positive entries of $\B$ sufficiently large and the negative entries sufficiently close to zero, which directly motivates the standardized matrix $\tilde{\A}$ in Corollary~\ref{cor:potential}.

\proof{Proof of Corollary~\ref{cor:potential}}
Suppose that $\HM(D) = 1$. By Theorem~\ref{thm:afr_ext}, $D$ satisfies GARP, so under $\A$, every cycle $\boldsymbol{c}$ in $[n]$ either consists solely of zero-weight transitions or contains at least one positive transition. In the first case, $\langle \mathbf{1}, \chi(\tilde{\A}, \boldsymbol{c}) \rangle = 0$. In the second case, $\boldsymbol{c}$ contains at most $n - 1$ negative transitions because it traverses at most $n$ distinct indices, and hence
\begin{equation*}
    \langle \mathbf{1}, \chi(\tilde{\A}, \boldsymbol{c}) \rangle \ge \frac{n}{n+1} - \frac{n-1}{n+1} = \frac{1}{n+1} > 0.
\end{equation*}
Thus every cycle in $[n]$ has nonnegative total weight under $\tilde{\A}$. By potential feasibility on cycles (see~\cite{gallai1958maximum}), there exists $\bphi \in \R^n$ satisfying $\phi_j - \phi_i \le \tilde{a}_{ij}$ for all $i, j \in [n]$.

Conversely, any such $\bphi$, together with $\blambda = \mathbf{1}$, is feasible for $\AI(\tilde{\A})$. Since $\sign(\tilde{\A}) = \sign(\A)$, Theorem~\ref{thm:afr_ext} yields $\HM(D) = 1$. \hfill \Halmos
\endproof

%%%%%%%%%%%%%%%%%%%%%%%%%%%%%%%%%%%%%%%%%%%%%%%%%
%%%%%%%%%%%% SECTION 4: Formulations %%%%%%%%%%%%
%%%%%%%%%%%%%%%%%%%%%%%%%%%%%%%%%%%%%%%%%%%%%%%%%
\section{MIP Formulation via Pairwise Convexification}\label{sec:formulations}

In this section, we leverage the potential-based characterization in Corollary~\ref{cor:potential} to derive a strong and compact formulation of $\HM(D)$. Specifically, in Section~\ref{sec:nonlinear}, we present a mixed 0-1 bilinear formulation of $\HM(D)$. In Section~\ref{sec:pairconv}, we derive the convex hull of each pair of bilinear inequalities. Finally, we present our MIP formulation in Section~\ref{sec:final}.
\subsection{A mixed 0-1 bilinear formulation}\label{sec:nonlinear}
To select a subset of the observations, we introduce binary variables $\x = (x_1, \ldots, x_n ) \in \{0,1\}^n$, where $x_i=1$ indicates that the observation $i$ is selected. We propose a mixed 0-1 bilinear formulation:
\begin{equation}\label{prog:nonlinear}
\begin{aligned}
    \max \quad & \frac{1}{n}\sum_{i=1}^n x_i \\
    \text{s.t.} \quad & \x\in\{0,1\}^n \quad \bphi\in[0,1]^n,\\
    & x_i(\phi_i-\phi_j+\tilde a_{ij})\ge0,\qquad i,j\in[n],
\end{aligned}
\end{equation}
where the bilinear constraint imposes the potential constraints $\phi_j-\phi_i\le\tilde a_{ij}$ for $j \in [n]$ whenever the observation $i$ is selected.

\begin{proposition}\label{prop:HM}
The formulation~\eqref{prog:nonlinear} solves $\HM(D)$:
\end{proposition}

\proof{Proof.}
For $\x\in\{0,1\}^n$, let $I_{\x}:=\{i\in[n]\mid x_i=1\}$ and $D_{\x}:=\{(\p^i,\q^i)\mid i\in I_{\x}\}$. For each $\x$ with $I_{\x}\ne\emptyset$, we show that $D_{\x}$ is rationalizable if and only if there exists a potential vector $\bphi\in[0,1]^n$ such that $(\x,\bphi)$ is feasible for~\eqref{prog:nonlinear}.

Suppose first that $D_{\x}$ is rationalizable. By Theorem~\ref{thm:afr_ext}, $D_{\x}$ satisfies GARP. Because $\sign(\tilde a_{ij})=\sign(a_{ij})$ for all $i,j\in I_{\x}$, every cycle $\boldsymbol{c}$ in $I_{\x}$ either has $\tilde a_{c_\ell, c_{\ell+1}}=0$ on every transition or contains at least one positive transition with weight $n/(n+1)$. Since a cycle in $I_{\x}$ has at most $|I_{\x}|\le n$ transitions and every negative transition has weight $-1/(n+1)$, any cycle containing a positive transition has total weight at least
\[
\frac{n}{n+1}-\frac{|I_{\x}|-1}{n+1}\ge\frac{1}{n+1}>0.
\]
Thus every cycle in $I_{\x}$ has nonnegative total length under the weight matrix $\tilde{\A}$, so by potential feasibility on cycles (see~\cite{gallai1958maximum}) there exist potentials $(\psi_i)_{i\in I_{\x}}$ satisfying $\psi_j-\psi_i\le\tilde a_{ij}$ for all $i,j\in I_{\x}$. In particular,
\[
\max_{k\in I_{\x}}\psi_k-\min_{k\in I_{\x}}\psi_k
=\max_{i,j\in I_{\x}}(\psi_j-\psi_i)\le\frac{n}{n+1}.
\]
Define $\bphi\in\R^n$ by
\begin{equation*}
    \phi_i:=
    \begin{cases}
        \psi_i-\max_{k\in I_{\x}}\psi_k+1, & i\in I_{\x},\\
        0, & i\notin I_{\x}.
    \end{cases}
\end{equation*}
Then $1/(n+1)\le\phi_i\le1$ for $i\in I_{\x}$ and $\phi_i=0$ for $i\notin I_{\x}$, so $\bphi\in[0,1]^n$. For $i,j\in I_{\x}$, we have $\phi_j-\phi_i=\psi_j-\psi_i\le\tilde a_{ij}$. For $i\in I_{\x}$ and $j\notin I_{\x}$, we have $\phi_j-\phi_i=-\phi_i\le-1/(n+1)\le\tilde a_{ij}$. Hence $x_i(\phi_i-\phi_j+\tilde a_{ij})\ge0$ holds for every $i\in I_{\x}$ and $j\in[n]$, while for $i\notin I_{\x}$ we have $x_i=0$ and the constraint holds trivially. Thus $(\x,\bphi)$ is feasible for~\eqref{prog:nonlinear}.

Conversely, suppose that $(\x,\bphi)$ is feasible for~\eqref{prog:nonlinear} with $I_{\x}\ne\emptyset$. For all $i,j\in I_{\x}$, setting $x_i=1$ in $x_i(\phi_i-\phi_j+\tilde a_{ij})\ge0$ yields $\phi_i-\phi_j+\tilde a_{ij}\ge0$. By Corollary~\ref{cor:potential},  $D_{\x}$ is rationalizable.   \Halmos
\endproof

\subsection{Pairwise convexification}\label{sec:pairconv}
One of the most effective techniques for deriving strong MIP formulations is convexification~\citep{conforti2014integer}. In our setting, we derive the convex hull of each pair of bilinear constraints. Namely, for $i < j$, consider the set
\begin{equation*}
    S_{ij} := \left\{(x_i,x_j,\phi_i,\phi_j) \in \{0,1\}^2 \times [0,1]^2 \ \middle|\ 
    \begin{aligned}
        x_i(\phi_i - \phi_j + \tilde a_{ij}) \ge 0, \quad x_j(\phi_j - \phi_i + \tilde a_{ji}) \ge 0
    \end{aligned}\right\}.
\end{equation*}

\begin{theorem}\label{thm:convex-hull}
For $i<j$, $\conv(S_{ij})$ consists of the points $(x_i,x_j,\phi_i,\phi_j)$ satisfying, 
\begin{subequations}\label{eq:paper-inequalities}
  \begin{align}
    &(x_i,x_j) \in [0,1]^2 \quad \text{ and } \quad x_i+x_j \leq 1   \for (i,j) \text{ s.t. } \tilde a_{ij}+\tilde a_{ji}<0 \label{eq:paper-inequalities-1} \\
    & \max\{-\tilde a_{ij},0\}x_i \le\phi_i\le1-\max\{-\tilde a_{ji},0\}x_j \label{eq:paper-inequalities-2}\\
    &\max\{-\tilde a_{ji},0\}x_j\le\phi_j\le1-\max\{-\tilde a_{ij},0\}x_i \label{eq:paper-inequalities-3}\\
    &(1-\tilde a_{ij})x_i-1\le\phi_i-\phi_j\le1-(1-\tilde a_{ji})x_j. \label{eq:paper-inequalities-4}
  \end{align}
\end{subequations}
\end{theorem}

\proof{Proof.}
Let $P_{ij}$ be the polytope defined by~\eqref{eq:paper-inequalities}. To show $\conv(S_{ij}) \subseteq P_{ij}$, it suffices by convexity of $P_{ij}$ to check that every $(x_i,x_j,\phi_i,\phi_j) \in S_{ij}$ lies in $P_{ij}$. Clearly $(x_i,x_j) \in [0,1]^2$. If $x_i = x_j = 1$, adding $\phi_i - \phi_j + \tilde a_{ij} \ge 0$ and $\phi_j - \phi_i + \tilde a_{ji} \ge 0$ gives $\tilde a_{ij} + \tilde a_{ji} \ge 0$; hence $\tilde a_{ij} + \tilde a_{ji} < 0$ forces $x_i + x_j \le 1$, proving~\eqref{eq:paper-inequalities-1}. Next, fix $(k,\ell) \in \{(i,j),(j,i)\}$. We verify
\begin{equation*}
    \max\{-\tilde a_{k\ell},0\}x_k \le \phi_k, \qquad \phi_\ell \le 1 - \max\{-\tilde a_{k\ell},0\}x_k, \qquad \phi_k - \phi_\ell \ge (1 - \tilde a_{k\ell})x_k - 1.
\end{equation*}
For $x_k = 0$, these reduce to $0 \le \phi_k$, $\phi_\ell \le 1$, and $\phi_k - \phi_\ell \ge -1$, which follow from $(\phi_i,\phi_j) \in [0,1]^2$. For $x_k = 1$, we have $\phi_k - \phi_\ell \ge -\tilde a_{k\ell} = (1 - \tilde a_{k\ell})x_k - 1$. Using $\phi_\ell \ge 0$ gives $\phi_k \ge -\tilde a_{k\ell}$ (so $\phi_k \ge \max\{-\tilde a_{k\ell},0\}$), and using $\phi_k \le 1$ gives $\phi_\ell \le 1 + \tilde a_{k\ell}$ (so $\phi_\ell \le \min\{1 + \tilde a_{k\ell},1\} = 1 - \max\{-\tilde a_{k\ell},0\}$). Applying these bounds for $(k,\ell) \in \{(i,j),(j,i)\}$ establishes~\eqref{eq:paper-inequalities-2}--\eqref{eq:paper-inequalities-4}.

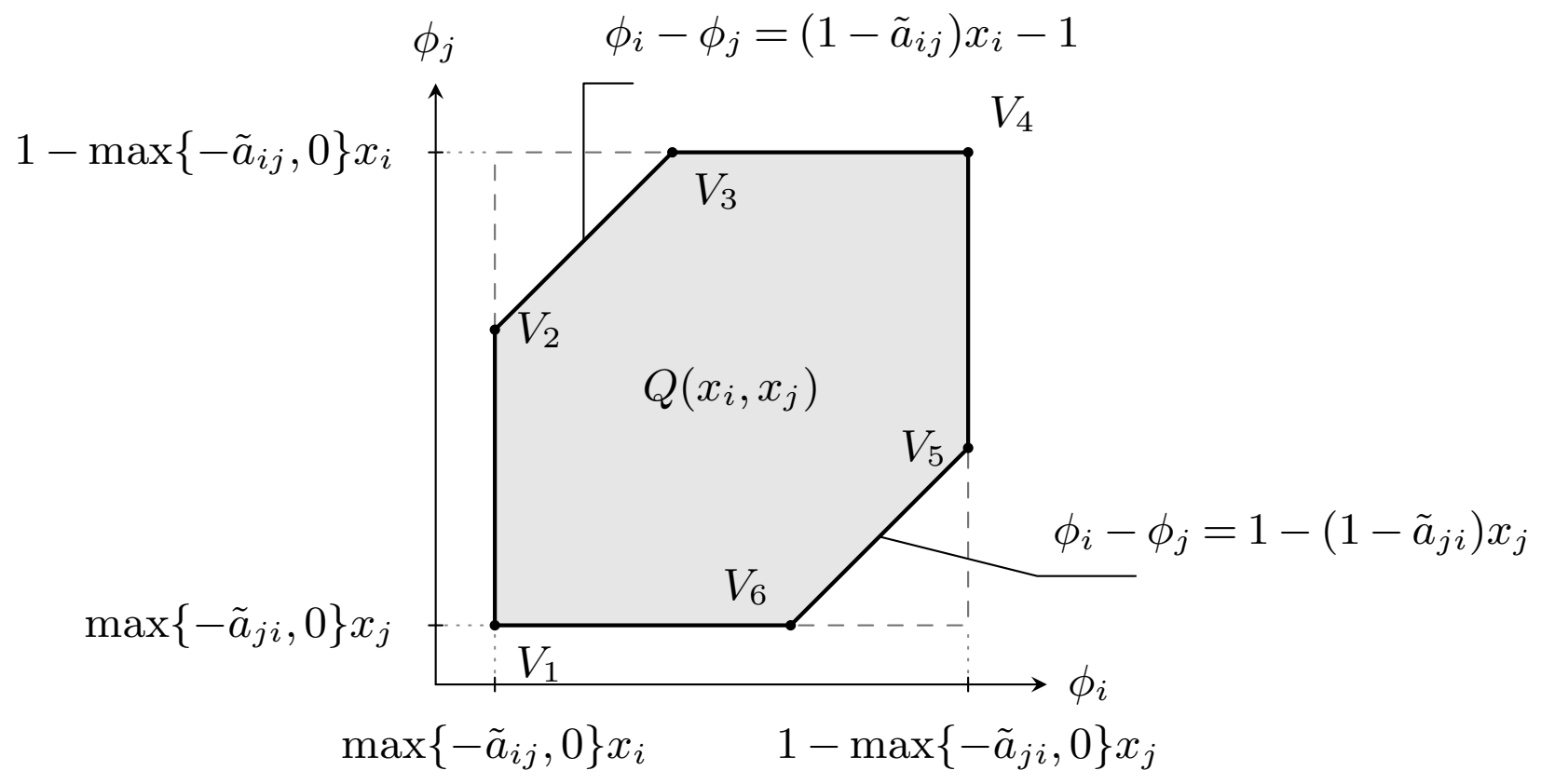
\begin{figure}[H]
    \centering
    \begin{tikzpicture}[scale=0.7,x=1cm,y=1cm,>=stealth,line cap=round,line join=round,every node/.style={font=\footnotesize}]
        \coordinate (V1) at (0.6,0.6);
        \coordinate (V2) at (0.6,3.6);
        \coordinate (V3) at (2.4,5.4);
        \coordinate (V4) at (5.4,5.4);
        \coordinate (V5) at (5.4,2.4);
        \coordinate (V6) at (3.6,0.6);

        \fill[black!10] (V1) -- (V2) -- (V3) -- (V4) -- (V5) -- (V6) -- cycle;
        \draw[dashed,black!55] (0.6,0.6) rectangle (5.4,5.4);
        \draw[thick] (V1) -- (V2) -- (V3) -- (V4) -- (V5) -- (V6) -- cycle;
        \node at (3,3) {$Q(x_i,x_j)$};

        \draw[->] (0,0) -- (6.2,0) node[right] {$\phi_i$};
        \draw[->] (0,0) -- (0,6.1) node[above] {$\phi_j$};
        \draw[dotted,black!45] (0.6,0) -- (V1) -- (0,0.6);
        \draw[dotted,black!45] (5.4,0) -- (5.4,0.6);
        \draw[dotted,black!45] (0,5.4) -- (0.6,5.4);
        \foreach \position in {0.6,5.4} {
            \draw (\position,0.07) -- (\position,-0.07);
            \draw (0.07,\position) -- (-0.07,\position);
        }
        \node[anchor=north] at (0.6,-0.2) {$\max\{-\tilde a_{ij},0\}x_i$};
        \node[anchor=north] at (5.4,-0.2) {$1-\max\{-\tilde a_{ji},0\}x_j$};
        \node[anchor=east] at (-0.2,0.6) {$\max\{-\tilde a_{ji},0\}x_j$};
        \node[anchor=east] at (-0.2,5.4) {$1-\max\{-\tilde a_{ij},0\}x_i$};

        \draw[thin] (1.5,4.5) -- (1.5,6.1) -- (2.0,6.1);
        \node[anchor=south west] at (1.5,6.15) {$\phi_i-\phi_j=(1-\tilde a_{ij})x_i-1$};
        \draw[thin] (4.5,1.5) -- (6.1,1.1) -- (7.1,1.1);
        \node[anchor=south west] at (6.05,1.1) {$\phi_i-\phi_j=1-(1-\tilde a_{ji})x_j$};

        \foreach \point in {V1,V2,V3,V4,V5,V6} {\fill (\point) circle (1.5pt);}
        \node[below right] at (V1) {$V_1$};
        \node[right] at (V2) {$V_2$};
        \node[below right] at (V3) {$V_3$};
        \node[above right] at (V4) {$V_4$};
        \node[left] at (V5) {$V_5$};
        \node[above left] at (V6) {$V_6$};
    \end{tikzpicture}
    \caption{The shaded set is $Q(x_i,x_j)$, shown with six distinct vertices.}
    \label{fig:pairwise-hull}
\end{figure}
We next prove the reverse inclusion $P_{ij} \subseteq \conv(S_{ij})$. Let $T$ be the set of points $(x_i,x_j)$ satisfying~\eqref{eq:paper-inequalities-1}. For a fixed $(x_i,x_j) \in T$, define $Q(x_i,x_j)$ as the set of points $(\phi_i,\phi_j)$ satisfying~\eqref{eq:paper-inequalities-2}, \eqref{eq:paper-inequalities-3}, and~\eqref{eq:paper-inequalities-4}. For every $(x_i,x_j) \in T$, $Q(x_i,x_j)$ is the hexagonal region depicted in Figure~\ref{fig:pairwise-hull}. To formalize its geometry, for $(k,\ell) \in \{(i,j),(j,i)\}$ we define the box boundaries
\begin{equation*}
    L_k(x_k) := \max\{-\tilde a_{k\ell}, 0\}x_k \quad \text{ and } \quad U_k(x_k) := 1 - \max\{-\tilde a_{ k \ell}, 0\}x_k,
\end{equation*}
and the diagonal cut intercepts
\begin{equation*}
    C(x_i) := (1-\tilde a_{ij})x_i - 1  \quad \text{ and } \quad D(x_j) := 1 - (1-\tilde a_{ji})x_j,
\end{equation*}
so that $Q(x_i,x_j)$ is the set of points $(\phi_i,\phi_j)$ satisfying $L_i(x_i) \le \phi_i \le U_j(x_j)$, $L_j(x_j) \le \phi_j \le U_i(x_i)$, and $C(x_i) \le \phi_i - \phi_j \le D(x_j)$. It can be verified that the six vertices, (allowing adjacent vertices to coincide on the boundary of $T$), can be  listed as follows:
\begin{align*}
    V_1(x_i,x_j) &:= \bigl(L_i(x_i),\, L_j(x_j)\bigr), &
    V_4(x_i,x_j) &:= \bigl(U_j(x_j),\, U_i(x_i)\bigr), \\
    V_2(x_i,x_j) &:= \bigl(L_i(x_i),\, L_i(x_i) - C(x_i)\bigr), &
    V_5(x_i,x_j) &:= \bigl(U_j(x_j),\, U_j(x_j) - D(x_j)\bigr), \\
    V_3(x_i,x_j) &:= \bigl(U_i(x_i) + C(x_i),\, U_i(x_i)\bigr), &
    V_6(x_i,x_j) &:= \bigl(L_j(x_j) + D(x_j),\, L_j(x_j)\bigr).
\end{align*}
Because each boundary line moves affinely with $(x_i,x_j)$, each mapping $V_k : T \to \R^2$ ($k \in [6]$) is affine on $T$.

Now, consider $(\bar{x}_i,\bar{x}_j, \bar{\phi}_i, \bar{\phi}_j) \in P_{ij}$, we express it as a convex combination of points in $S_{ij}$. Clearly, $(\bar{\phi}_i,\bar{\phi}_j) \in Q(\bar{x}_i, \bar{x}_j)$, and thus can be expressed a convex combination of $V_k(\bar{x}_i, \bar{x}_j)$ for $k \in [6]$, that is $(\bar{\phi}_i,\bar{\phi}_j) = \sum_{k\in [6]} \alpha_kV_k(\bar{x}_i, \bar{x}_j)$, where $\alpha_k \ge 0$ for all $k \in [6]$ and $\sum_{k \in [6]} \alpha_k = 1$.
Moreover, letting $\operatorname{vert}(T) \subseteq \{0,1\}^2$ denote the set of vertices of $T$, we have $(\bar{x}_i, \bar{x}_j) \in T$, so $(\bar{x}_i, \bar{x}_j)$ can be expressed as a convex combination of $\operatorname{vert}(T)$: $    (\bar{x}_i, \bar{x}_j) = \sum_{\z \in \operatorname{vert}(T)} \beta_{\z} \z$, where $\beta_{\z} \ge 0$ for all $\z = (z_i,z_j) \in \operatorname{vert}(T)$ and $\sum_{\z \in \operatorname{vert}(T)} \beta_{\z} = 1$. Define the product weights $\gamma_{k,\z} := \alpha_k \beta_{\z} \ge 0$ for $(k,\z) \in [6] \times \operatorname{vert}(T)$, which satisfy $\sum_{k \in [6]} \sum_{\z \in \operatorname{vert}(T)} \gamma_{k,\z} = 1$. Using the affinity of each $V_k$ on $T$, we compute
\begin{align*}
    \sum_{k \in [6]} \sum_{\z \in \operatorname{vert}(T)} \gamma_{k,\z} \bigl(z_i, z_j, V_k(z_i, z_j)\bigr)
    &= \sum_{k \in [6]} \alpha_k \Biggl( \sum_{\z \in \operatorname{vert}(T)} \beta_{\z} \z,\; \sum_{\z \in \operatorname{vert}(T)} \beta_{\z} V_k(z_i, z_j) \Biggr) \\
    &= \sum_{k \in [6]} \alpha_k \bigl(\bar{x}_i, \bar{x}_j, V_k(\bar{x}_i, \bar{x}_j)\bigr) 
     = (\bar{x}_i, \bar{x}_j, \bar{\phi}_i, \bar{\phi}_j).
\end{align*}
It can be verified that $\bigl(z_i, z_j, V_k(z_i, z_j)\bigr) \in S_{ij}$ for every $k \in [6]$. Hence, we obtain $(\bar{x}_i, \bar{x}_j, \bar{\phi}_i, \bar{\phi}_j) \in \conv(S_{ij})$ and completes the proof. \hfill \Halmos

\subsection{Final formulation}\label{sec:final}
Now, we are ready to present our formulation:
\begin{equation*}\label{eq:compact}
    \max \Biggl\{ \frac{1}{n}\sum_{i=1}^n x_i \Biggm| \x\in\{0,1\}^n,\ (x_i,x_j,\phi_i,\phi_j) \text{ satisfies}~\eqref{eq:paper-inequalities} \for   i<j \Biggr\}.
\tag{\textsc{Compact}}
\end{equation*}
As we show next,~\eqref{eq:compact} is an MIP formulation with $n$ binary variables, $n$ continuous variables, and $O(n^2)$ linear inequalities.

\begin{corollary}\label{cor:compact}
\eqref{eq:compact} is an MIP formulation of $\HM(D)$.
\end{corollary}
\proof{Proof.}
By Proposition~\ref{prop:HM}, $\HM(D)$ is equivalent to $\max \bigl\{ \frac{1}{n}\sum_{i=1}^n x_i \bigm| (x_i,x_j,\phi_i,\phi_j) \in S_{ij} \text{ for all } i < j \bigr\}$. Observe that for $i < j$, $S_{ij} = \bigl\{(x_i,x_j,\phi_i,\phi_j) \in \conv(S_{ij}) \bigm| (x_i,x_j) \in \{0,1\}^2 \bigr\}$. By Theorem~\ref{thm:convex-hull}, replacing each $(x_i,x_j,\phi_i,\phi_j) \in S_{ij}$ with $(x_i,x_j) \in \{0,1\}^2$ and~\eqref{eq:paper-inequalities} yields an exact formulation of $\HM(D)$. \hfill \Halmos
\endproof

%%%%%%%%%%%%%%%%%%%%%%%%%%%%%%%%%%%%%%%%%%%%%%%
%%%%%%%%%%%% SECTION 5: Extensions %%%%%%%%%%%%
%%%%%%%%%%%%%%%%%%%%%%%%%%%%%%%%%%%%%%%%%%%%%%%
% \section{Extensions}\label{sec:extensions}

% In this section, we extend our combinatorial characterization of rationalizability and the compact HM formulation to three settings: general budget sets, finite choice sets, and symmetric rationalizability.  The first two extend the structure of the observed datasets, while the third considers a different notion of rationalizability. The key ingredient underlying these extensions is Theorem~\ref{thm:afr_ext}, which allows us to transform existing Afriat-type characterizations into combinatorial Afriat-type inequalities and subsequently derive compact HM formulations. These results demonstrate the broader applicability of our combinatorial approach beyond the classical rationalizability setting.

%%%%%%%%%%%%%%%%%%%%%%%%%%%%%%%%%%%%%%%%%%%%%%%%%%%%%%%%%%
%%%%%%%%%%%% SECTION 6: Computational Results %%%%%%%%%%%%
%%%%%%%%%%%%%%%%%%%%%%%%%%%%%%%%%%%%%%%%%%%%%%%%%%%%%%%%%%
\section{Computational Results}\label{sec:computational-results}
We benchmark \eqref{eq:compact} against the state-of-the-art formulation of \citet{demuynck2023computing} (DR), detailed in Appendix~\ref{app:dr-formulation}. All experiments were conducted in Julia~1.12.6 using JuMP~1.31.2 and Gurobi~13.0.3 on a Windows~11 machine with an Intel Core i7-12700H processor and 16\,GB of RAM. Each run was restricted to a single thread, a 3,600-second MIP time limit, and an 8\,GB soft memory limit.

\paragraph{Test instances.}
Each instance comprises $n$ observations on $m$ goods, with prices $p^i_j$ drawn independently from $\operatorname{Unif}[1,2]$ and total expenditure normalized to one in every observation. We generate positive weights $Y^i_j$ independently of prices and set quantities via $q^i_j = Y^i_j / \bigl(p^i_j \sum_{k=1}^m Y^i_k\bigr)$ under two data-generating processes:
\begin{itemize}
    \item \emph{Perturbed Cobb--Douglas choices:} Adapting the two-good CES weights of \citet[Appendix~III]{fisman2007individual}, we set unperturbed expenditure shares $w_1 = 3/4$ and $w_j = 1/[4(m-1)]$ for $j > 1$. Following \citet[Appendix~C]{blundell2003nonparametric}, we perturb these shares by drawing $Y^i_j \sim \operatorname{Beta}\bigl(\kappa w_j, \kappa(1-w_j)\bigr)$ independently across goods and observations, where $\kappa = (m+1)/\lambda - 1$. Larger values of $\lambda$ increase the variance of the draws around $w_j$, producing stronger deviations from exact Cobb--Douglas rationality.
    \item \emph{Purely random choices:} We draw $Y^i_j \sim \operatorname{Exp}(1)$ independently across goods and observations. The normalized weights are uniformly distributed on the unit simplex, inducing the uniform distribution on the budget hyperplane studied by \citet{bronars1987power}.
\end{itemize}
We generate ten random instances for each parameter combination $(n, m, \lambda)$ in the Cobb--Douglas suite and each pair $(n, m)$ in the purely random suite.

\paragraph{Performance metrics.}
Let $z^{\mathrm{best}}$ denote the cardinality of the largest rationalizable subset found by either formulation; then $z^{\mathrm{best}}/n$ equals $\HM(D)$ whenever either formulation certifies optimality and provides a lower bound otherwise. Reported MIP time measures Gurobi's branch-and-bound solve time, with runs hitting the time or memory limit recorded at 3,600 seconds. We define the LP relaxation gap as $100(z^{\mathrm{LP}} - z^{\mathrm{best}})/z^{\mathrm{best}}$, where $z^{\mathrm{LP}}$ is the optimal value of the root LP relaxation, and the final optimality gap as $100(U - L)/L$ for $L > 0$, where $U$ and $L$ are the upper bound and incumbent value at termination. Overall means are weighted by the number of instances contributing to each group average.

\paragraph{Summary of results.}
Tables~\ref{tab:cd-summary} and~\ref{tab:pure-summary} summarize the computational results. Across all 179 reported instances, \eqref{eq:compact} certifies optimality on 159 instances compared with 101 for DR, and strictly dominates DR in both LP relaxation tightness and solve time:
\begin{itemize}
    \item \emph{Perturbed Cobb--Douglas choices (Table~\ref{tab:cd-summary}):} \eqref{eq:compact} solves all 119 instances to optimality with an overall mean MIP time of $33.04$\,s, achieving a $42.5$-fold speedup over DR ($1{,}404.98$\,s). Its pairwise convexification yields a nearly exact LP relaxation, reducing the overall mean LP gap from $6.27\%$ to $0.12\%$. For $n = 1{,}000$ and $\lambda = 0.020$, \eqref{eq:compact} reduces the mean LP gap from $20.67\%$ to $0.08\%$ ($m = 2$) and from $14.02\%$ to $0.80\%$ ($m = 5$), solving all 40 instances with $n = 1{,}000$ within three minutes on average, whereas every DR run at $n = 1{,}000$ hits the time limit, leaving large optimality gaps.
    \item \emph{Purely random choices (Table~\ref{tab:pure-summary}):} \eqref{eq:compact} tightens the LP relaxation in every parameter group, reducing the overall mean LP gap from $229.76\%$ to $69.18\%$, and solves 40 of 60 instances (including all instances with $m = 2$ up to $n = 1{,}000$) compared with 22 of 60 for DR. On the harder $m = 5$ instances with $n \in \{500, 1{,}000\}$ where neither formulation proves optimality within one hour, \eqref{eq:compact} reduces the mean final gap from $254.63\%$ to $4.71\%$ at $n = 500$ and from $50{,}200.00\%$ to $111.91\%$ at $n = 1{,}000$.
\end{itemize}

\begin{table}[!htbp]
\centering
\small
\renewcommand{\arraystretch}{0.92}
\caption{Perturbed Cobb--Douglas choices.}
\label{tab:cd-summary}
\setlength{\tabcolsep}{2.5pt}
\begin{tabular*}{0.96\linewidth}{@{\extracolsep{\fill}}rrrrrrrrrr@{}}
\toprule
$n$ & $m$ & $\lambda$ & Mean HM & \multicolumn{2}{c}{Mean MIP time (s)} & \multicolumn{2}{c}{Mean LP gap (\%)} & \multicolumn{2}{c}{Mean final gap (\%)} \\
\cmidrule(lr){5-6}\cmidrule(lr){7-8}\cmidrule(lr){9-10}
 & & & & DR & Compact & DR & Compact & DR & Compact \\
\midrule
200 & 2 & 0.005 & 0.9890 & 2.4647 & \textbf{1.5669} & 1.1152 & \textbf{0.0000} & 0.0000 & 0.0000 \\
200 & 2 & 0.020 & 0.9465 & 2.9825 & \textbf{1.2549} & 5.6794 & \textbf{0.0000} & 0.0000 & 0.0000 \\
200 & 5 & 0.005 & 0.9980 & 1.6317 & \textbf{1.1188} & 0.2015 & \textbf{0.0000} & 0.0000 & 0.0000 \\
200 & 5 & 0.020 & 0.9640 & 2.5155 & \textbf{1.1874} & 3.7518 & \textbf{0.0000} & 0.0000 & 0.0000 \\
\midrule
500 & 2 & 0.005 & 0.9778 & 596.4449 & \textbf{8.0627} & 2.2751 & \textbf{0.0000} & 0.0000 & 0.0000 \\
500 & 2 & 0.020 & 0.8882 & 665.6474 & \textbf{9.4014} & 12.6011 & \textbf{0.0563} & 0.0000 & 0.0000 \\
500 & 5 & 0.005 & 0.9960 & 421.5737 & \textbf{8.1319} & 0.4022 & \textbf{0.0000} & 0.0000 & 0.0000 \\
500 & 5 & 0.020 & 0.9252 & 685.6574 & \textbf{20.1399} & 8.0889 & \textbf{0.4300} & 0.0000 & 0.0000 \\
\midrule
1000 & 2 & 0.005 & 0.9513 & 3600 & \textbf{57.2016} & 5.1271 & \textbf{0.0000} & 31566.6667 & \textbf{0.0000} \\
1000 & 2 & 0.020 & 0.8289 & 3600 & \textbf{63.4304} & 20.6688 & \textbf{0.0848} & 34066.6667 & \textbf{0.0000} \\
1000 & 5 & 0.005 & 0.9911 & 3600 & \textbf{40.0808} & 0.8993 & \textbf{0.0101} & 30733.3333 & \textbf{0.0000} \\
1000 & 5 & 0.020 & 0.8771 & 3600 & \textbf{182.4325} & 14.0195 & \textbf{0.8024} & 31066.6667 & \textbf{0.0000} \\
\midrule
\multicolumn{3}{l}{Overall} & 0.9441 & 1404.9809 & \textbf{33.0423} & 6.2691 & \textbf{0.1163} & 10708.6835 & \textbf{0.0000} \\
\bottomrule
\end{tabular*}
\par\smallskip
\begin{minipage}{0.96\linewidth}
\footnotesize
Each group contains ten instances, except for $(n,m,\lambda)=(500,2,0.005)$, for which one instance was excluded from both formulations' averages. The DR separation parameter was $\delta\approx5.96\times10^{-7}$, below Gurobi's feasibility tolerance of $10^{-6}$. DR returned 487 observations at a reported zero gap, but the subset violated GARP. Compact certified an optimum of 486. This illustrates the numerical sensitivity of DR when $\delta$ is small.
\end{minipage}

\medskip
\caption{Purely random choices.}
\label{tab:pure-summary}
\setlength{\tabcolsep}{2.5pt}
\begin{tabular*}{0.96\linewidth}{@{\extracolsep{\fill}}rrrrrrrrrrr@{}}
\toprule
$n$ & $m$ & \multicolumn{1}{c}{\multirow{2}{*}{\shortstack{Mean\\[0pt]$z^{\mathrm{best}}/n$}}} & \multicolumn{2}{c}{Solved} & \multicolumn{2}{c}{Mean MIP time (s)} & \multicolumn{2}{c}{Mean LP gap (\%)} & \multicolumn{2}{c}{Mean final gap (\%)} \\
\cmidrule(lr){4-5}\cmidrule(lr){6-7}\cmidrule(lr){8-9}\cmidrule(lr){10-11}
 & & & DR & Compact & DR & Compact & DR & Compact & DR & Compact \\
\midrule
200 & 2 & 0.5225 & 10 & 10 & 79.1659 & \textbf{1.3985} & 91.6457 & \textbf{4.3699} & 0.0000 & 0.0000 \\
200 & 5 & 0.4125 & 8 & 10 & 1863.8459 & \textbf{16.0716} & 142.9803 & \textbf{27.0518} & 0.8889 & \textbf{0.0000} \\
\midrule
500 & 2 & 0.4456 & 4 & 10 & 3246.9765 & \textbf{26.3217} & 124.6701 & \textbf{15.9279} & 89.6224 & \textbf{0.0000} \\
500 & 5 & 0.3210 & 0 & 0 & 3600 & 3600 & 212.2337 & \textbf{59.0747} & 254.6281 & \textbf{4.7108} \\
\midrule
1000 & 2 & 0.3838 & 0 & 10 & 3600 & \textbf{292.3353} & 160.8674 & \textbf{32.4493} & 40813.8889\textsuperscript{*} & \textbf{0.0000} \\
1000 & 5 & 0.1430 & 0 & 0 & 3600 & 3600 & 646.1401 & \textbf{276.2198} & 50200.0000\textsuperscript{*} & \textbf{111.9080} \\
\midrule
\multicolumn{2}{l}{Overall} & 0.3714 & 22 & 40 & 2664.9980 & \textbf{1256.0212} & 229.7562 & \textbf{69.1822} & 7265.3068\textsuperscript{*} & \textbf{19.4365} \\
\bottomrule
\end{tabular*}
\par\smallskip
\begin{minipage}{0.96\linewidth}
\footnotesize
Each group contains ten instances. All 20 DR runs with $n=1{,}000$ reached the 8\,GB soft memory limit before one hour and are recorded at 3,600 seconds each.
\par\noindent\textsuperscript{*}\,The mean excludes runs with a zero incumbent, leaving 6 of 10 runs for $(n,m)=(1{,}000,2)$, 2 of 10 for $(1{,}000,5)$, and 48 of 60 overall.
\end{minipage}
\end{table}

\section{Conclusion}\label{sec:conclusion}
In this paper, we revisited the Afriat inequalities through a combinatorial lens, establishing a scale-free Afriat characterization. Leveraging this characterization, we formulated the computation of $\HM(D)$ as a mixed 0-1 bilinear program and derived the exact convex hull of each pairwise bilinear subsystem. The resulting compact MIP formulation requires only $O(n)$ variables and $O(n^2)$ linear inequalities, eliminates the big-$M$ numerical sensitivities of existing formulations, and achieves substantial improvements in LP relaxation tightness and solve times on both perturbed Cobb--Douglas and purely random datasets.

\FloatBarrier
\bibliographystyle{informs2014} % outcomment this and next line in Case 1
\bibliography{refs} 

% the ecompanion here
% Options are (1) APPENDIX (with or without general title) or 
%             (2) APPENDICES (if it has more than one unrelated sections)
% Outcomment the appropriate case if necessary
%
\clearpage
% \begin{APPENDIX}{title}
% \end{APPENDIX}
%
%   or 
%

\clearpage

\ECSwitch

\ECDisclaimer
%%%%%%%%%%%%%%%%%%%%%%%%%%%%%%%%%%%%%%%%%%%%%%%%%%%%%%%%%%

%%% Main head for the e-companion
\ECHead{Additional Discussions, Missing Proofs, Formulations, and Numerical Studies}

 \section{The State-of-the-Art Formulation}\label{app:dr-formulation}

\citet{demuynck2023computing} formulate the HMI using the following characterization of GARP. A consumption dataset $D$ with revealed preference matrix $\A$ satisfies GARP if and only if there exists a set of numbers $\{\phi_i \}_{i \in [n]} \subset [0,1]$ such that for all $i,j \in [n]$,
\begin{equation*}
    \begin{aligned}
        a_{ij} \le 0 \quad & \text{implies} \quad \phi_i \ge \phi_j, \\
        a_{ij} < 0 \quad & \text{implies} \quad \phi_i > \phi_j.
    \end{aligned}
\end{equation*}

Their formulation is generally constructed in three steps. First, they constructed $n^2$ binary variables $\{ U_{ij} \}_{i,j \in [n]}$ as indicators for the condition $\phi_i \ge \phi_j$, while introducing $\epsilon \in (0,1/n)$ to handle $\phi_i > \phi_j$. Specifically, the condition ``$U_{ij} = 1$ if and only if $\phi_i \ge \phi_j$'' is ensured by the following inequalities:
\begin{align*}
    \phi_i - \phi_j & \le -\epsilon + 2U_{ij}, \\
    U_{ij} - 1 & \le \phi_i - \phi_j.
\end{align*}

Second, they introduced fixed parameters $\alpha$ and $\delta$ such that $\alpha > \max_{i \in [n]} \{\langle \p^i, \q^i \rangle \}$ and $0< \delta < \min \left\{ \min_{i,j \in [n]} \{ \langle \p^i, \q^j \rangle \}, \min_{i,j \in [n]} \{a_{ij} \mid a_{ij} > 0 \}  \right\}$. Here we assume that all observed bundles are nonzero and use the convention $\min\varnothing=+\infty$. With the help of these parameters, the following inequalities can capture the condition ``$a_{ij} \le 0$ implies $U_{ij} = 1$ and $a_{ij} < 0$ implies $U_{ji} = 0$'':
\begin{align*}
    - a_{ij} \le -\delta + \alpha U_{ij}, \\
    \alpha(U_{ji} - 1) \le a_{ij}.
\end{align*}

Finally, they formulated $\HM(D)$ by introducing binary variables $x_i$ to indicate whether the $i$-th observation is selected. Their complete formulation is presented as follows:

\begin{equation*}
\begin{aligned}
      \max \quad & \frac{\sum_{i=1}^n x_i}{n} \\
    \text{s.t.} \quad & \x \in \{0,1\}^n,\ \bphi \in [0,1]^n \text{ and } \boldsymbol{U} \in \{0,1\}^{n \times n} \\
    & \phi_i - \phi_j \le -\epsilon + 2U_{ij} \quad \for i,j \in [n] \\
    & U_{ij} - 1 \le \phi_i - \phi_j \quad \for i,j \in [n] \\ 
    & -a_{ij} \le -\delta + \alpha [U_{ij} + (1- x_i)] \quad \for i,j \in [n] \\ 
    & \alpha(U_{ji} + x_i -2) \le a_{ij} \quad \for i, j \in [n]
\end{aligned}
\end{equation*}

Since this formulation contains $O(n^2)$ binary variables, it becomes computationally expensive for standard solvers to find an optimal solution when $n$ is large.

 \section{Extensions}\label{sec:extensions}

In this section, we extend our combinatorial characterization of rationalizability and the compact HM formulation to three settings: general budget sets, finite choice sets, and symmetric rationalizability.  The first two extend the structure of the observed datasets, while the third considers a different notion of rationalizability. The key ingredient underlying these extensions is Theorem~\ref{thm:afr_ext}, which allows us to transform existing Afriat-type characterizations into combinatorial Afriat-type inequalities and subsequently derive compact HM formulations. These results demonstrate the broader applicability of our combinatorial approach beyond the classical rationalizability setting.

\subsection{General Budget Sets}

Following \citet{forges2009afriat}, consider observations $D=\{(\q^i,B_i)\}_{i=1}^n$, where $\q^i$ is chosen from the budget set $B_i=\{\q\in\R_{\ge0}^m:g_i(\q)\le0\}$. Each $g_i:\R_{\ge0}^m\to\R$ is a continuous and increasing function, with $g_i(\q^i)=0$. Here, increasing means that the function value strictly increases whenever all coordinates increase. Since $g_i$ is not necessarily linear or convex, the budget set $B_i$ may be nonlinear and nonconvex. The dataset $D$ is rationalizable if there exists a continuous and increasing utility function $u:\R_{\ge0}^m\to\R$ such that
\begin{equation*}
    \q^i\in\arg\max_{\q\in B_i}u(\q),\qquad i\in[n].
\end{equation*}
Let $\HM(D)$ denote the maximum fraction of observations in $D$ that can be rationalized by such a utility function.

By \citet[Proposition~3 and its proof]{forges2009afriat}, $D$ is rationalizable if and only if there exist $\bphi\in\R^n$ and $\blambda\in\R_{>0}^n$ satisfying
\begin{equation*}
    \phi_i-\phi_j+\lambda_i g_i(\q^j)\ge0,\qquad i,j\in[n].
\end{equation*}
Thus, setting $a_{ij}:=g_i(\q^j)$ gives the system $\AI(\A)$, with $a_{ii}=g_i(\q^i)=0$. The same characterization holds for every selected subdataset. Using the standardization in Corollary~\ref{cor:potential}, Proposition~\ref{prop:HM} and Corollary~\ref{cor:compact} therefore give an exact formulation of $\HM(D)$ through \eqref{eq:compact}. 

\subsection{Finite Choice Sets}

Following \citet{cosaert2015revealed}, consider a dataset $D=\{(\q^i,B_i)\}_{i=1}^n$. Each $B_i=\{\boldsymbol{b}_i^1,\ldots,\boldsymbol{b}_i^{K_i}\}\subseteq\R_{\ge0}^m$ is a finite budget set with $K_i$ elements, and $\q^i\in B_i$ is the observed choice. A utility function $u:\R_{\ge0}^m\to\R$ is weakly monotone if $u(\q)\ge u(\q')$ whenever $\q\ge\q'$ coordinatewise, with strict inequality whenever every coordinate of $\q$ exceeds that of $\q'$. Here, $D$ is rationalizable if there exists a continuous, weakly monotone utility function $u$ such that $\q^i\in\arg\max_{\q\in B_i}u(\q)$ for every $i\in[n]$.

For this setting, \citet[Section~3.1]{cosaert2015revealed} define the revealed preference matrix $\A=(a_{ij})$ by
\begin{equation*}
    a_{ij}:=\min_{k\in[K_i]}\max_{\ell\in[m]}\bigl(q_\ell^j-b_{i\ell}^k\bigr),\qquad i,j\in[n].
\end{equation*}
They show that $D$ is rationalizable if and only if there exist $\bphi\in\R^n$ and $\blambda\in\R_{>0}^n$ satisfying
\begin{equation*}
    \phi_i-\phi_j+\lambda_i a_{ij}\ge0,\qquad i,j\in[n].
\end{equation*}
These inequalities give precisely $\AI(\A)$, so the same matrix arguments apply to $D$ and every selected subdataset.

Since $\q^i\in B_i$, we have $a_{ii}\le0$. If $a_{ii}<0$, some bundle in $B_i$ strictly exceeds $\q^i$ in every coordinate, so observation $i$ cannot belong to any rationalizable subdataset. To include diagonal comparisons in \eqref{eq:compact}, we impose the lower bound on $\phi_i-\phi_j$ in \eqref{eq:paper-inequalities-4} also for $j=i$. This gives $(1-\tilde a_{ii})x_i\le1$, which forces $x_i=0$ when $a_{ii}<0$ because $x_i$ is binary. If $a_{ii}=0$, the inequality reduces to the existing bound $x_i\le1$. With these diagonal inequalities, \eqref{eq:compact} computes $\HM(D)$, the maximum fraction of observations in a rationalizable subdataset, using the standardization in Corollary~\ref{cor:potential}.

\subsection{Symmetric Rationalizability}\label{sec:symmetric}

We next consider rationalization by a symmetric utility function, following \citet{chambers2018note}. Let $\Sigma_m$ be the set of all permutations $\sigma:[m]\to[m]$. For $\boldsymbol{v}\in\R^m$, the permuted vector $\sigma\boldsymbol{v}$ moves the coordinate in position $\ell$ to position $\sigma(\ell)$, so $(\sigma\boldsymbol{v})_{\sigma(\ell)}=v_\ell$. A utility function $u:\R_{\ge0}^m\to\R$ is symmetric if $u(\sigma\q)=u(\q)$ for every $\q\in\R_{\ge0}^m$ and $\sigma\in\Sigma_m$. We say that $D=\{(\p^i,\q^i)\}_{i=1}^n$ is symmetrically rationalizable if there exists a continuous, concave, strictly increasing, and symmetric utility function $u$ such that
\begin{equation*}
    \q^i\in\arg\max\left\{u(\q)\;\middle|\;\q\in\R_{\ge0}^m,\ \langle\p^i,\q\rangle\le\langle\p^i,\q^i\rangle\right\},\qquad i\in[n].
\end{equation*}
Let $\HM_{\mathrm{sym}}(D)$ be the maximum fraction of observations contained in a symmetrically rationalizable subdataset.

Proposition~1 of \citet{chambers2018note} states that $D$ is symmetrically rationalizable if and only if there exist $\bphi\in\R^n$ and $\blambda\in\R_{>0}^n$ satisfying
\begin{equation}\label{eq:symmetric-afriat-expanded}
    \phi_j-\phi_i\le\lambda_i\langle\eta\p^i,\sigma\q^j-\eta\q^i\rangle,\qquad i,j\in[n],\ \sigma,\eta\in\Sigma_m.
\end{equation}
These inequalities use one utility level $\phi_i$ and one positive multiplier $\lambda_i$ per original observation, shared across all its permutations.

We next reduce the permutation comparisons in \eqref{eq:symmetric-afriat-expanded} to one coefficient for each pair $(i,j)$. Since simultaneous permutations preserve inner products,
\begin{align*}
    \langle\eta\p^i,\sigma\q^j-\eta\q^i\rangle
    &=\langle\eta\p^i,\sigma\q^j\rangle-\langle\eta\p^i,\eta\q^i\rangle\\
    &=\langle\p^i,\eta^{-1}(\sigma\q^j)\rangle-\langle\p^i,\q^i\rangle.
\end{align*}
As $\sigma$ and $\eta$ vary, $\eta^{-1}(\sigma\q^j)$ ranges over all permutations of $\q^j$. Because $\lambda_i>0$, the inequalities in \eqref{eq:symmetric-afriat-expanded} for a fixed pair $(i,j)$ are equivalent to the inequality with the smallest right-hand side. We therefore define the revealed preference matrix $\A^{\mathrm{sym}}=(a^{\mathrm{sym}}_{ij})$ by
\begin{equation}\label{eq:symmetric-coefficients}
\begin{aligned}
    a^{\mathrm{sym}}_{ij}&:=\min_{\tau\in\Sigma_m}\langle\p^i,\tau\q^j\rangle-\langle\p^i,\q^i\rangle\\
    &=\langle\p^{i,\uparrow},\q^{j,\downarrow}\rangle-\langle\p^i,\q^i\rangle,\qquad i,j\in[n],
\end{aligned}
\end{equation}
where $\p^{i,\uparrow}$ lists the price coordinates in nondecreasing order and $\q^{j,\downarrow}$ lists the quantity coordinates in nonincreasing order. The second equality follows from the rearrangement inequality, so the minimum can be computed by sorting. Consequently, $D$ is symmetrically rationalizable if and only if there exist $\bphi\in\R^n$ and $\blambda\in\R_{>0}^n$ satisfying
\begin{equation}\label{eq:symmetric-afriat}
    \phi_i-\phi_j+\lambda_i a^{\mathrm{sym}}_{ij}\ge0,\qquad i,j\in[n].
\end{equation}
This is the system $\AI(\A^{\mathrm{sym}})$, so the potential characterization and the selection arguments apply with the matrix $\A^{\mathrm{sym}}$.

We finally explain why negative diagonal entries do not prevent the use of our formulation. The identity permutation gives $a^{\mathrm{sym}}_{ii}\le0$. Setting $j=i$ in \eqref{eq:symmetric-afriat} gives $\lambda_i a^{\mathrm{sym}}_{ii}\ge0$, so $a^{\mathrm{sym}}_{ii}<0$ excludes observation $i$ from every symmetrically rationalizable subdataset. Using the standardization in Corollary~\ref{cor:potential}, we include diagonal comparisons by imposing the lower bound on $\phi_i-\phi_j$ in \eqref{eq:paper-inequalities-4} also for $j=i$. This gives $(1-\tilde a^{\mathrm{sym}}_{ii})x_i\le1$, which forces $x_i=0$ when $a^{\mathrm{sym}}_{ii}<0$ because $x_i$ is binary. If $a^{\mathrm{sym}}_{ii}=0$, the inequality reduces to the existing bound $x_i\le1$. Thus, \eqref{eq:compact}, with these diagonal inequalities included, computes $\HM_{\mathrm{sym}}(D)$ exactly.

The formulation uses $n$ binary selection variables, $n$ continuous potentials, and $O(n^2)$ explicitly specified linear inequalities. Together with the sorting computation of $\A^{\mathrm{sym}}$, this gives an exact formulation of polynomial size without enumerating permutations or cycles.

\end{document}